\documentclass[12pt]{amsart}

\newtheorem{lemma}{Lemma}[section]
\newtheorem{theorem}[lemma]{Theorem}
\newtheorem{corollary}[lemma]{Corollary}

\newtheorem{definition}[lemma]{Definition}
\newtheorem{condition}[lemma]{Condition}

\newcommand{\rr}{\mathbb{R}}
\newcommand{\pp}{\mathbb{P}}
\newcommand{\cc}{\mathbb{C}}

\newcommand{\Oo}{\mathcal{O}}
\newcommand{\Rr}{\mathcal{R}}
\newcommand{\grad}{\mbox{{\bf grad}}\, }
\newcommand{\hull}{\mbox{{\bf hull}}\, }

\author[C. Simpson]{Carlos Simpson}
\address{CNRS, Laboratoire J. A. Dieudonn\'e\\ Universit\'e de Nice-Sophia Antipolis\\
06108 Nice, Cedex 2, France}
\email{carlos@math.unice.fr}
\urladdr{http://math.unice.fr/$\sim$carlos/} 
\title[Asymptotics for connections]{Asymptotics for general connections at infinity}

\begin{document}


\begin{abstract}
For a standard path of connections going to a generic point at infinity in the
moduli space $M_{DR}$ of connections on a compact Riemann surface, we show that the Laplace transform of the
family of monodromy matrices has an analytic continuation with locally finite branching. In particular the 
convex subset representing the exponential growth rate of the monodromy is a polygon, whose vertices are in 
a subset of points described explicitly in terms of the spectral curve. Unfortunately we don't get any information about
the size of the singularities of the Laplace transform, which is why we can't get asymptotic expansions
for the monodromy. 
\end{abstract}




\maketitle


\section{Introduction}
\label{introduction}

We study the asymptotic behavior of the monodromy of connections
near a general point at $\infty$ in the space $M_{DR}$ of connections on a compact Riemann surface $X$.
We will consider a path of connections of the form $(E,\nabla + t\theta )$ which approaches
the boundary divisor transversally at the  point on the boundary of $M_{DR}$ corresponding to a
general Higgs bundle $(E,\theta )$. 
By some meromorphic gauge transformations in \S \ref{pullback} we reduce to the case of a family of
connections of the form $d+B+tA$. This is very similar to what was treated in \cite{abmade} except that 
here our matrix $B$ may have poles. 
We import the vast majority of our techniques directly from there. 
The difficulty posed by the poles of $B$ is the new phenomenon which is treated here.
We are not able to get results as good as the precise asymptotic expansions of \cite{abmade}.
We just show in Theorem \ref{main} (p.\ \pageref{main}) that if $m(t)$ denotes the family of monodromy
or transport matrices for a given path, then the Laplace transform $f(\zeta )$ of $m$ has an analytic
continuation with locally finite singularities over the complex plane (see Definition \ref{lfb},
p.\ \pageref{lfb}). The singularities are what
determine the asymptotic behavior of $m(t)$. The upside of this situation is that since we are aiming for less, 
we can considerably simplify certain parts of the argument. 
What we don't know is the behavior of $f(\zeta )$ near
the singularities: the main question left open is whether $f$ has polynomial growth at the
singularities, and if so, to what extent the generalized Laurent series can be calculated from
the individual terms in our integral expression for $f$. 

We can get some information about where the singularities are.  Fix a general point $(E,\theta )$.
Recall from \cite{Hitchin} \cite{Hitchin2} \cite{Donagi} \cite{Kanev} 
that the {\em spectral curve} $V$ is the subset of points in $T^{\ast}(X)$ corresponding to eigenforms
of $\theta$.  We have a proper mapping $\pi : V\rightarrow X$. In the case of a general point, $V$ is smooth and
the mapping has only simple ramification points. Also there is a tautological one-form 
$$
\alpha \in 
H^0(V, \pi ^{\ast}\Omega ^1_X)\subset H^0(V,\Omega ^1_V).
$$
Finally there is a line bundle $L$ over $V$ such that $E\cong \pi _{\ast}(L)$ and $\theta$ corresponds to the
action of $\alpha$ on the direct image bundle. This is all just a geometric version of the diagonalization of
$\theta$ considered as a matrix over the function field of $X$. 

Let $\Rr \subset X$ denote the subset of points over which the spectral curve is ramified, that is 
the image of the set of branch points of $\pi$. It is the set of {\em turning points} of our singular perturbation
problem. Suppose $p$ and $q$ are points in $X$ joined by a path $\gamma$.
A {\em piecewise homotopy lifting} of $\gamma$ to the spectral curve $V$ consists of a collection of paths
$$
\tilde{\gamma}=
\{ \tilde{\gamma}_i\} _{i=1,\ldots , k}
$$ 
such that each $\tilde{\gamma}_i$ is a continuous path in $V$,
and such that if we denote by $\gamma _i:= \pi \circ \tilde{\gamma}_i$ the image paths in $V$, then
$\gamma _1$ starts at $p$, $\gamma _k$ ends at $q$, and for $i=1,\ldots , k-1$, the endpoint of
$\gamma _i$ is equal to the starting point of $\gamma _{i+1}$ and this is a point in $\Rr$.
Among these there is a much more natural class of paths which are the {\em continuous homotopy liftings}, namely those
where the starting point of $\tilde{\gamma} _i$ is equal to the starting point of $\tilde{\gamma}_{i+1}$ (which is
not necessarily the case for a general piecewise lifting).

Denote by $\Sigma (\gamma )\subset \cc$ the set of integrals of the tautological form $\alpha$ 
along piecewise homotopy liftings of $\gamma$, i.e. the set of complex numbers of the form 
$$
\sigma = \int _{\tilde{\gamma}}\alpha := \sum _{i=1}^k \int _{\tilde{\gamma}_i}\alpha .
$$
Let $\Sigma ^{\rm cont}(\gamma )$ be the subset of integrals along the continuous homotopy liftings.
The following is the statement of Theorem \ref{main} augmented with a little bit of information about where
the singularities are.

\begin{theorem}
\label{mainplus}
Let $p,q$ be two points on $X$, and let $\gamma$ denote a path from $p$ to $q$. Let $\{ (E, \nabla + t\theta )\}$
denote a curve of connections cutting the divisor $P_{DR}$ at a general point $(E,\theta )$ and let
$(V,\alpha , L)$ denote the spectral data for this Higgs bundle. Let $m(t)$ be the function (with values in
$Hom (E_p,E_q)$) whose value at $t\in \cc$ is the transport matrix for the connection $\nabla + t\theta $
from $p$ to $q$ along the path $\gamma$. Let $f(\zeta )$ denote the Laplace transform of $m$. Then 
$f$ has an analytic continuation with locally finite singularities over the complex plane. The set of
singularities which are ever encountered is a subset of the set $\Sigma (\gamma )\subset \cc$
of integrals of the tautological form along piecewise homotopy liftings defined above. 
\end{theorem}

It would have been much nicer to be able to say that the set of singularities is contained in
$\Sigma ^{\rm cont}(\gamma )$, however I don't see that this is necessarily the case. However, it might
be that the singularities in $\Sigma ^{\rm cont}(\gamma )$ have a special form different from the others.
This is an interesting question for further research.

Ths first singularities which are encountered in the analytic continuation of $f$
determine the growth rate of $m(t)$ in a way which we briefly formalize. 
Suppose that
$m(t)$ is an entire function with exponentially bounded growth. 
We say that $m(t)$ is {\em rapidly decreasing in a sector}, if for some (open) 
sector of complex numbers going to $\infty$, there is $\varepsilon >0$ giving
a bound of the form $|m(t)|\leq e^{-\varepsilon |t|}$. 
Define the {\em hull} of $m$ by
$$
\hull(m) := \{ \zeta \in \cc \, \mbox{s.t.}\,  e^{-\zeta t} m(t) \mbox{not rapidly decreasing in any sector} \}.
$$
It is clear from the definition that the set of
$\zeta$ such that $e^{-\zeta t} m(t)$ is rapidly decreasing in some sector, is open. Therefore $\hull (m)$ 
is closed. It is also not too hard to see that it is convex (see \S \ref{proof}). Note that
the hull is defined entirely in terms of the growth rate of the function $m$. 

\begin{corollary}
\label{hullcor}
In the situation of Theorem \ref{mainplus}, the hull of $m$ is a finite convex polygon with at least two vertices,
and all of its vertices  are contained in $\Sigma (\gamma )$. 
\end{corollary}

The above results fall into the realm of {\em singular perturbation theory} for systems of ordinary
differential equations, which goes back at least to Liouville. A steady stream of progress in this theory
has led to a vast literature which we don't attempt completely to cover here (and which the reader can
explore by using internet and database search techniques, starting for example from the authors mentionned
in the bibliography).  

Recall that following \cite{BalianBloch}, Voros and Ecalle
looked at these questions from the viewpoint of ``resurgent functions'' \cite{Voros} \cite{Voros2} \cite{VorosBourbaki}
\cite{Ecalle1} \cite{EcalleOrsay} \cite{EcalleAbordable} \cite{Boutet} \cite{CNP} \cite{DDP}. 
In the terminology  of Ecalle's article in \cite{Boutet}, the singular perturbation problem we are considering here is
an example of {\em co-equational resurgence}. 
Our approach is very related to this viewpoint, though self-contained. We use a notion of analytic continuation
of the Laplace transform \ref{lfb} which is a sort of weak version of resurgence, like that 
used in \cite{DDP} and \cite{CNP}. The elements of our 
expansion \ref{formula} are what Ecalle calls the ``elementary monomials'' and the trees which appear in \S \ref{description}
are related to 
{\em (co)moulds (co)arborescents}, see \cite{Boutet}. Conversion properties related to the trees
have been discussed in \cite{EcalleVallet} (which is on the subject of KAM theory \cite{GentileMastropietro}). 
The relationship with integrals on a spectral curve was explicit in \cite{DDthesis}, \cite{DDP}. 
The works \cite{VorosBourbaki}, Ecalle's article in \cite{Boutet}, and \cite{DDP},
raise a number of questions about how to prove resurgence for certain classes of singular perturbation problems
notably some arising in quantum mechanics. 
A number of subsequent articles treat these questions; I haven't been able to include everything here but some
examples are  \cite{EcalleVallet}, \cite{DDP2}, \cite{DelabaereHowls}, \ldots (and apparently \cite{ZinnJustin}).
In particular \cite{DelabaereHowls} discuss extensively the way in which the singularities of the Laplace
transform determine the asymptotic behavior of the original function, specially in the case of the
kinds of integrals which appear as terms in the decomposition \ref{formula}. 

There are a number of other currents of thought about the problem of singular perturbations. 
It is undoubtedly important to pursue the relationship with all of these. 
For example, the study initiated in \cite{BerkNevinsRoberts} and continuing with several articles in \cite{Boutet},
as well as the more modern \cite{AokiKawaiTakei} (also Prof. Kawai's talk at this conference) 
indicates that there is an intricate and fascinating geometry in the propagation of the Stokes phenomenon. 
And on the other hand it would be good to understand the relationship with the local study of
turning points such as in \cite{RamisEtAl}, \cite{Stenger}. The article \cite{DDP2} incorporates some
aspects of all of these approaches, and one can see \cite{BenderEtAl} for a physical perspective. 
Also works on Painlev\'e's equations and isomonodromy such as \cite{CostinCostin} 
\cite{InabaIwasakiSaito} \cite{SanguinettiWoodhouse} \cite{Woodhouse} are probably relevant.

Even though he doesn't appear in the references of \cite{abmade}, the ideas of J.-P. Ramis indirectly
had a profound influence on that work (and hence on the present note). 
This can be traced to at least two inputs as follows: 
\newline
(1) I had previously followed G. Laumon's course
about $\ell$-adic Fourier transform, which was partly inspired by the corresponding notions in complex function 
theory, a subject in which Ramis (and Ecalle, Voros, \ldots ) had a great influence; and
\newline
(2) at the time of writing \cite{abmade} I was following N. Katz's course about exponential
sums, where again much of the inspiration came from Ramis' work (which Katz mentionned very often) on irregular 
singularities.  

Thus I would like to take this opportunity to thank Jean-Pierre for inspiring such a rich mathematical context. 

I would also like to thank the several participants who made interesting remarks and posed interesting questions.
In particular F. Pham pointed out that
it would be a good idea to look at what the formula for the location of the singularities actually said,
leading to the statement of Theorem \ref{main} in its above form. I haven't been able to treat other
suggestions (D. Sauzin, \ldots ), such as looking at the differential equation satisfied by $f(\zeta )$.

\section{The compactified moduli space of connections}
\label{compactified}

Let $X$ be a smooth projective curve over the complex numbers $\cc$. Fix $r$ and suppose $E$ is a 
vector bundle of rank $r$ over $X$. A {\em connection} (by which we mean an algebraic one) on $E$
is a $\cc$-linear morphism of sheaves $\nabla : E \rightarrow E\otimes _{\Oo}\Omega ^1_X$ satisfying
the Leibniz rule $\nabla (ae) = (da) e + a\nabla (e)$. If $p$ and $q$ are points joined by a path
$\gamma : [0,1]\rightarrow X$, $\gamma (0)=P$, $\gamma (1)=Q$ then local solutions of $\nabla (e)=0$ continue
along $\gamma$, giving a {\em transport matrix} $m_{\gamma}(E,\nabla ) : E_P\rightarrow E_Q$. 
The transport matrix, our main object of study, is the fundamental solution of a linear system of ODE's. 
If $E$ is a trivial bundle
(which will always be the case at least on a Zariski open subset of $X$ containing $\gamma$) then
there is a formula for the transport matrix as a sum of iterated integrals \cite{Chen} \cite{Hain}. A modified
version of this formula is basic to the argument below, although we mostly refer to \cite{abmade} for the details
of that part of the argument.

Recall that we have
a moduli space $M_{DR}$ of rank $r$ vector bundles with integrable connection on $X$ \cite{moduli},
which has a
compactification $M_{DR}\subset \overline{M}_{DR}$ constructed as follows. 
A {\em Higgs bundle} is a  pair $(E,\theta )$ where $\theta : E \rightarrow E\otimes _{\Oo}\Omega ^1_X$
is an $\Oo _X$-linear bundle map (rather than a connection) \cite{Hitchin} \cite{Hitchin2} \cite{hbls}, which 
is {\em semistable of degree $0$} if $E$ has degree zero and if any sub-Higgs bundle has degree $\leq 0$.
In fact for any $\lambda \in \mathbb{A}^1_{\cc}$ we can look at the notion of {\em vector bundle with 
$\lambda$-connection} \cite{DeligneLet}---related in an obvious way to the notion of 
singular perturbation---which is a pair $(E,\nabla )$ of a bundle plus a connection-like operator satisfying
Leibniz' rule with a factor of $\lambda$ in front of the first term. For $\lambda =0$ this is just a Higgs bundle
and for any $\lambda \neq 0$ the operator $\lambda ^{-1}\nabla$ is a connection.

With these definitions, there is a moduli space \cite{naht} \cite{hfnac} \cite{moduli}
$M_{Hod}\rightarrow {\bf A}^1$ for vector bundles with $\lambda$-connection, $\lambda \in {\bf A}^1$.
The fiber over $\lambda  =0$ is the moduli space $M_{Dol}$ for semistable Higgs bundles of degree 
zero, whereas for any $\lambda \neq 0$ the fiber is isomorphic to $M_{DR}$. 

The Higgs-bundle moduli space has a subvariety $M_{Dol}^{\rm nil}$ parametrizing the Higgs bundles $(E,\theta )$ such that
$\theta$ is nilpotent as an $\Omega ^1_X$-valued endomorphism of $E$. Let $M_{Dol}^{\ast}$ denote
the complement of $M_{Dol}^{\rm nil}$ in $M_{Dol}$ and let  
$M_{Hod}^{\ast}$ denote
the complement of $M_{Dol}^{\rm nil}$ in $M_{Hod}$. Then the algebraic group ${\bf G}_m$ acts on
$M_{Hod}$ preserving all of the above subvarieties, and the compactification is obtained as the quotient
\cite{naht} \cite{hfnac}
$$
\overline{M}_{DR} := M_{Hod}^{\ast} /{\bf G}_m.
$$
The complement of $M_{Dol}$ in $M_{Hod}$ (which is also 
the complement of $M_{Dol}^{\ast}$ in $M_{Hod}^{\ast}$) is isomorphic to 
$M_{DR}\times {\bf G}_m$ and this gives the embedding $M_{DR}\hookrightarrow \overline{M}_{DR}$.
The complementary divisor is given by 
$$
P_{DR} = M_{Dol}^{\ast} /{\bf G}_m.
$$
In conclusion, this means that the points at $\infty$ in $\overline{M}_{DR}$ correspond to 
equivalence classes of semistable, degree $0$, non-nilpotent Higgs bundles $(E, \theta )$ 
under the equivalence relation
$$
(E,\theta )\cong (E, u\theta )
$$
for any $u\in {\bf G}_m$.

Recall that the moduli space $M_{Dol}$ is an irreducible algebraic variety \cite{moduli}, so $P_{DR}$ is also
irreducible. The general point therefore corresponds to a ``general'' Higgs bundle $(E,\theta )$ (in what follows
we often forget to add the adjectives ``semistable, degree $0$'').
For a general point, the {\em spectral curve} of $\theta$ (described in more detail in the section
after next) is an irreducible curve with ramified map to $X$,
such that the ramification points are all of the simplest type. 

We should note that Arinkin \cite{Arinkin} \cite{ArinkinLysenko} has
defined a finer compactification by modifying the notion of $\lambda$-connection,
and this is taken up by Inaba, Iwasaki and Saito \cite{InabaIwasakiSaito}.

\section{Curves going to infinity}
\label{curves}

The moduli spaces considered above are coarse only. In an etale neighborhood of the generic point, though,
they are fine and smooth. At a general point of the divisor $P_{DR}$,
both  $\overline{M}_{DR}$ and $P_{DR}$ are smooth. Thus we can look for a curve cutting $P_{DR}$ transversally 
at a general point. Such a curve may be obtained by taking the projection of a curve in 
$M_{Hod}$ cutting $M_{Dol}$ at a general point. In turn, this amounts to giving a family 
$(E_c,\nabla _c)$ where $\nabla _c$ is a $\lambda (c)$-connection, parametrized by $c\in C$
for some curve $C$. In an etale neighborhood of the point $\lambda =0$, the function
$\lambda (c)$ should be etale. Note also that $(E_0,\nabla _0)$ should be a 
general semistable Higgs bundle of degree zero. 

The easiest way to obtain such a curve is as follows: let $(E,\theta )$ be a general Higgs bundle, stable
of degree zero. The bundle $E$ is stable as a vector bundle (since stability is an open condition and it
certainly holds on the subset of Higgs bundles with $\theta =0$, so it holds at general points). In particular
$E$ supports a connection $\nabla$ and we can set 
$$
\nabla _{\lambda} := \lambda \nabla + \theta
$$
for $\lambda\in {\bf A}^1$. Here the parameter is $\lambda$ itself.
The subset ${\bf G}_m\subset {\bf A}^1$ corresponds to points which are mapped into 
$M_{DR}$, and indeed the vector bundle with connection corresponding to the above $\lambda$-connection
is 
$$
(E, \nabla + t\theta ) \;\;\; , \; t= \lambda ^{-1}.
$$
The map actually extends to a map from ${\bf A}^1$ into $M_{DR}$ for the other coordinate chart 
${\bf A}^1$  providing a neighborhood at $\infty$ in ${\bf P}^1$. 
In conclusion, the family of connections $\{ (E, \nabla + t\theta )  \}$ corresponds to a morphism
$$
{\bf P}^1 \rightarrow \overline{M}_{DR}
$$
sending $t\in {\bf A}^1$ into $M_{DR}$, sending the point $t=\infty$ to a general point in the divisor
$P_{DR}$, and the curve is transverse to the divisor at that point. This type of curve was called a {\em pencil of
connections} by Losev and Manin \cite{LosevManin}. 

We will look only at curves of the above form. It should be possible to obtain similar results for 
other curves cutting $P_{DR}$ transversally at a general point, but that is left as a problem for future study. 

We will investigate the asymptotic behavior of the monodromy representations of the connections
$(E, \nabla + t\theta )$ as $t\rightarrow \infty$. Recall that the {\em Betti moduli space}
$M_B$ is the moduli space for representations of $\pi _1(X)$ up to conjugation, and we have an
analytic isomorphism $M_{DR}^{\rm an}\cong M_B^{\rm an}$ sending a connection to its monodromy
representation. We will look at the asymptotics of the resulting analytic curve 
${\bf A}^1\rightarrow M_B$. 

In order to set things up it will be useful to fix a basepoint $p\in X$ and a trivialization $\tau : E_p\cong
\cc ^r$. Then for any $\gamma \in \pi _1(X,x)$ we obtain the monodromy matrix
$$
\rho (E, \nabla + t\theta ,\tau , \gamma )\in GL(r,\cc ).
$$
Of course the monodromy matrices don't directly give functions on the moduli space $M_B$ of representations,
because they depend on the choice of trivialization $\tau$.
However, one has the Procesi coordinates (see Culler and Shalen \cite{CullerShalen} and Procesi \cite{Procesi}) 
which are certain polynomials in the monodromy matrices (for
several $\gamma$ at once) which are invariant under change of trivialization and give an embedding 
of the Betti moduli space $M_B$ into an affine space. We will look at the asymptotic
behavior of the monodromy matrices, but the resulting asymptotic information will also hold for
any polynomials (see Corollary \ref{polynomials}), 
and in particular for the Procesi coordinates. This will give asymptotic information
about the image curve in $M_B$. 

Notationally it is easier to start right out considering the transport matrices between points $p$ and $q$.
In any case, the functions we shall consider, be they the matrix coefficients of the monodromy $\rho$
or some other polynomials in these or the transport matrices, 
will be entire functions $m(t)$ on the complex line $t\in \cc$. 
We will be looking to characterize their asymptotic properties. 

The method we will use is the same as the method already used in \cite{abmade} to treat exactly
this question, for a more special class of curves going to infinity in $M_{DR}$. In that book
was treated the case of families of connections $(E,\nabla + t\theta )$ where
$$
E=\Oo _X^r,\;\; \nabla = d + B, \;\; \theta = A
$$
with $A$ and $B$ being $r\times r$ matrices of one-forms on $X$ such that $A$ is diagonal and
$B$ contains only zeros on the diagonal. In \cite{abmade}, a fairly precise description of the
asymptotic behavior of the monodromy was obtained. It was also indicated how one should be able to
reduce to this case in general; we shall explain that below. The only problem is that in the
course of this reduction, one obtains the special situation but with $B$ being a matrix of one-forms
which has some poles on $X$. In this case the exact method used in \cite{abmade} breaks down.

The purpose of the present paper is to try to remedy this situation as far as possible. We change very 
slightly the method (essentially by taking the more canonical gradient flows of the functions $\Re g_{ij}$ rather
than the flows defined in Chapter 3 of \cite{abmade}, and also stopping the flows before arriving at the
poles of $B$). However, we don't obtain the full results of \cite{abmade}, namely we can show an analytic continuation
result for the Laplace transform of $m(t)$ (this Laplace transform is explained in more detail below),
however we don't get good bounds or information about the singularities of the Laplace transform other
than that they are locally finite sets of points. In particular we obtain information about the growth
rate of $m(t)$ but not asymptotic expansions. 

Even in order to obtain the analytic continuation, a much more detailed examination of the dynamics
generated by the general method of \cite{abmade} is necessary. This is the main body of the present
paper (see Theorem \ref{mainbound}). For the remainder of the technique we mostly refer to \cite{abmade}. 

Thus while we treat a much more general type of curve going to infinity than was treated
in \cite{abmade}, we obtain a weaker set of results for these curves.  This leaves open the
difficult question of what kinds of singularities the Laplace transforms have, and thus what
type of asymptotic expansion we can get for $m(t)$. 

\section{Genericity results for the spectral data}
\label{genericity}

Before beginning to look more closely at the monodromy representations, we will consider some
properties of general points $(E,\theta )$ on $P_{DR}$, best expressed
in terms of the {\em spectral curve} \cite{Hitchin} \cite{Hitchin2} \cite{Donagi} \cite{Kanev}
\cite{DDthesis} \cite{SanguinettiWoodhouse}.

Suppose $(E,\theta )$ is a Higgs bundle. Suppose $P\in X$ and $v\in T_PX$; 
then we obtain the fiber $E_P$ which is a vector space
of rank $r$, with an endomorphism $\theta _P(v)\in End (E_P)$. We say that
$P$ is {\em singular} if $\theta_P(v)$ has an eigenvalue (i.e. zero of the characteristic polynomial)
of multiplicity $\geq 2$.
It is more natural to look at the {\em eigenforms} of $\theta$ obtained by dividing out the vector $v$.
The eigenforms are elements of the cotangent space $T^{\ast}_PX = (\Omega ^1_X)_P$. 

We say that a singular point $P$ is {\em generic} if there is exactly one eigenform
of multiplicity $\geq 2$; if it has multiplicity exactly $2$; and if the
two eigenforms $\alpha ^{\pm}$ of $\theta $ which come together at $P$, may be expressed 
in a neighborhood with coordinate $z$ as
$$
\alpha ^{\pm} = c dz \pm az^{1/2}dz + \ldots .
$$
The condition that all singular points are generic is a Zariski open condition on the moduli space of Higgs bundles.

Suppose $P$ is a generic singular point. The eigenforms give a set of $r-1$ distinct elements of $T^{\ast}_PX$,
consisting of  the values of the multiplicity-one eigenvalues of $\theta$ at $P$, plus the leading 
term $cdz$ for the pair $\alpha ^{\pm}$. Call this set
$EF _P$.  We say that $P$ is {\em non-parallel} if 
$EF _P$, viewed as a subset of the
real two-dimensional space $T^{\ast}_PX$, doesn't have any colinear triples, nor
any quadruples of points defining two parallel lines.

In terms of a coordinate $z$ at $P$ we can write the elements of $EF _P$ as
$$
\alpha _i(P) = a_i dz
$$
with $a_i$ being distinct complex numbers, and say $a_1=c$ in the previous formulation. 
Then $P$ is non-parallel if and only if the set of
$a_i \in \cc \cong \rr ^2$ doesn't have any colinear triples or parallel quadruples.
In turn this is equivalent to saying that the angular coordinates of the complex  numbers
$a_i-a_j$ are distinct. 

\begin{lemma}
\label{nonparallelgen}
The set of Higgs bundles $(E,\theta )$ such that the singularities are generic and satisfy the
non-parallel condition, is a dense real Zariski-open subset of the moduli space. 
\end{lemma} 
\begin{proof}
The condition of being non-parallel is a real Zariski open condition. In particular, the condition 
that all singular points be generic and non-colinear, holds in the complement
of a closed real algebraic subset of the moduli space. Therefore, if there is one such point then 
the set of such points is a dense real Zariski open subset. 

To show that there is one point $(E,\theta )$ such that all of the singular points are generic and
non-parallel, we can restrict to the case where $E= \Oo ^{\oplus r}$ is a trivial bundle. 
In this case, $\theta$ corresponds to a matrix of holomorphic one-forms on $X$. We will consider a matrix
of the form $A + \lambda B$ with $A$ diagonal having entries $\alpha _i$, and $B$ is off-diagonal with 
$\lambda$ small.  The singular points are perturbations of the points where $\alpha _i(P)= \alpha _j(P)$.
A simple calculation with a $2\times 2$ matrix shows that the singularities are generic in this case.
In order to obtain the non-colinear condition, it suffices to have that for a point $P$ where
$\alpha _i(P)=\alpha _j(P)$, the subset of $r-1$ valuse of all the $\alpha _k(P)$ is non-parallel. 

For a general choice of the $\alpha _k$, this is the case. Suppose we are at a point $P$ where
$\alpha _1(P)=\alpha _2(P)$ for example. Then moving the remaining $\alpha _k$ for $k\geq 3$
shows that the remaining points are general with respect to the first one. A set of $r-1$ points
such that the last $r-2$ are general with respect to the first one (whatever it is), satisfies
the non-parallel condition.
\end{proof}

\begin{lemma}
\label{lb}
If $(E,\theta )$ is generic in the sense of the previous lemma, then the spectral curve $V$ is
actually an irreducible smooth curve sitting in the cotangent bundle $T^{\ast}X$. There is a line bundle
$L$ on $V$ such that $E\cong \pi _{\ast} (L)$ and $\theta$ is given by multiplication by the
tautological one-form over $V$. 
\end{lemma}
\begin{proof}
The genericity condition on the way the eigenforms come together at any point where the multiplicity
is $\geq 2$, guarantees that at any point where the projection $\pi : V\rightarrow X$ is not
locally  etale, the curve $V$ is a smooth ramified covering of order $2$ in
the usual standard form. This shows that $V$ is smooth. It is irreducible, 
because this is so for at least some points (for example the deformations used in the
previous proof) and Zariski's connectedness implies that in a connected family of smooth projective
curves if one is irreducible then all are. For connectedness of the family we use the
irreducibility of the moduli space of Higgs bundles cf \cite{moduli}. The last statement is 
standard in theory of spectral curves \cite{Hitchin} \cite{Hitchin2} \cite{Donagi} \cite{Kanev}.
\end{proof}

{\em Remark:} Once $p$ and $q$ are fixed, then for general $\theta$ the endpoints $p,q$ will not be contained in the set
$\Rr$ of turning points. 

\section{Pullback to a ramified covering and gauge transformations}
\label{pullback}

Fix a general Higgs bundle $(E,\theta )$ on $X$. By taking a Galois completion of the spectral curve of
$\theta$ and Galois-completing a 
further two-fold ramified covering if necessary, we can obtain a ramified Galois
covering 
$$
\varphi  : Y\rightarrow X
$$
such that the pullback Higgs field $\varphi ^{\ast}$ has a full set of eigen-one-forms defined on $Y$;
and such that the ramification powers over singular points of $\theta$ are divisible by $4$. 

We have one-forms $\alpha _1,\ldots , \alpha _r$ and line sub-bundles 
$$
L_1,\ldots , L_r \subset \varphi ^{\ast}E
$$
such that at a general point of $Y$ we have
$$
\psi : \varphi ^{\ast}E\cong L_1\oplus \ldots \oplus L_r
$$
with $\varphi ^{\ast}\theta$ represented by the diagonal matrix with entries $\alpha _i$. Note that 
$\varphi ^{\ast}\theta$ preserves $L_i$ (acting there by multiplication by $\alpha _i$)
globally on $Y$. However, the isomorphism $\psi$ will only be meromorphic, and also the $L_i$ are
of degree $<0$. Choose modifications $L'_i$ of $L_i$ (see Lemma \ref{modif} below, also the modifications
are made only over singular points) such that $L'_i$ is 
of degree zero, and set
$$
E':= L'_1\oplus \ldots \oplus L'_r.
$$
Let $\theta '$ denote the diagonal Higgs field with entries $\alpha _i$ on $E'$. 
Let $\nabla '$ be a diagonal flat connection on $E'$. 
We have a meromorphic map
$$
\psi : E \rightarrow E',
$$
and 
$$
\psi \circ \varphi ^{\ast}\theta \circ \psi ^{-1} = \theta '.
$$

Suppose  now that $\nabla$ was a connection on $E$, giving a connection $\varphi ^{\ast}\nabla$ on 
$\varphi ^{\ast}E$. We can write
$$
\psi \circ \varphi ^{\ast} \nabla \circ \psi ^{-1} = \nabla ' + \beta 
$$
with $\beta$ a meromorphic  section of $End (E') \otimes _{\Oo} \Omega ^1_Y$. 

A transport matrix of $(E,\nabla + t\theta )$ may be recovered as a transport matrix for
the pullback bundle on $Y$. Indeed if $\gamma$ is a path in $X$ going from $p$ to $q$ then
it lifts to a path going from a lift $p'$ of $p$ to a lift $q'$ of $q$.  Thus it suffices to 
look at the problem of the asymptotics for transport matrices for the pullback family
$$
\{ (\varphi ^{\ast}E,\varphi ^{\ast}\nabla + t\varphi ^{\ast}\theta ) \} .
$$
We may assume that $p$ and $q$  are not singular points of $\theta$, so $p'$ and $q'$ will not
be singular points of $\varphi ^{\ast}\theta $. Then the transport matrices for this family
are conjugate (by a conjugation which is constant in $t$) to the transport matrices for the
family 
$$
\{ (E', \nabla ' + \beta + t\theta ' )\} .
$$

\begin{lemma}
\label{modif}
In the above situation, the modifications $L'_i$ of $L_i$ may be chosen so that the diagonal entries of
$\beta$ are  holomorphic. Furthermore the poles of the remaining entries of $\beta$ are restricted to the
points lying over singular points in $X$ for the original Higgs field $\theta$ (the ``turning points''). 
\end{lemma}
\begin{proof}
Note first that, by definition, away from the singular points of $\theta$ the eigen-one-forms are distinct
so the eigenvectors form a basis for $E$, in other words the direct sum decomposition $\psi$ is an
isomorphism at these points. Thus $\psi$ only has poles over the singular points of $\theta$ (hence the
same for $\beta$).

We will describe a choice of $L'_i$ locally at a singular point. 

Look now in a neighborhood of a point $P'\in Y$, lying over a singular point $P\in X$. Let 
$z'$ denote a local coordinate at $P'$ on $Y$, with $z$ a local coordinate at $P$ on $X$
and with 
$$
z = (z')^m.
$$
Our assumption on $Y$ was that $m$ is divisible by $4$. In fact we may as well assume that $m=4$
since raising to a further power doesn't modify the argument. Thus we can write
$$
z' = z^{1/4}.
$$
There are two eigenforms of $\theta$ which come together at $P$. Suppose that their lifts are $\alpha _1$ and
$\alpha _2$. Then near $P'$ we can write 
$$
\varphi ^{\ast}E = U \oplus L_3 \oplus \ldots \oplus L_r
$$
where $U$ is the rank two subbundle of $\varphi ^{\ast}E$ corresponding to eigenvalues 
$\alpha _1$ and $\alpha _2$. The direct sum decomposition  is holomorphic at $P'$ because
the other eigenvalues of $\theta$ were distinct at $P$ and different from the two singular ones
(of course after the pullback all of the eigenforms have a value of zero at $P'$ but the decomposition 
still holds nonetheless).

Now we use a little bit more detailed information about spectral curves for Higgs bundles: the general 
$(E,\theta )$ is obtained as the direct image of a line bundle on the spectral curve (Lemma \ref{lb}). This means that
locally near $P$ there is a two-fold branched covering with coordinate $u=z^{1/2}$ such that the 
rank $2$ subbundle of $E$ corresponding to the singular values looks like the direct
image of the trivial bundle on the covering, and the $2\times 2$ piece of $\theta$ looks like 
the action of multiplication by $udz = 2u^2du$. The direct image, considered as a module over the
series in $z$, is just the series in $u$. One can obtain a basis by looking at the odd and even 
powers of $u$: the basis vectors are $e_1=1$ and $e_2=u$. In these terms we have
$$
\theta e_1 = e_2dz; \;\;\; \theta e_2 = z e_1dz.
$$
Thus the $2\times 2$ singular part of $\theta$ has matrix
$$
\left( \begin{array}{cc}
0& z \\
1 & 0 \end{array} \right) dz.
$$
Pulling back now to the covering $Y$ which is locally $4$-fold, we have a basis for $U$ in which
$$
\varphi ^{\ast}\theta |_U = 
\left( 
\begin{array}{cc}
0& (z')^7 \\
(z')^3 & 0 \end{array} \right) dz'.
$$
On the other hand, since up until now our decomposition is holomorphic, the pullback connection
$\varphi ^{\ast}\nabla$ may be written (in terms of our basis for $U$ plus trivializations of the $L_i$
for $i\geq 3$) as $d+B'$ where $B'$ is a holomorphic matrix of one-forms. Since the basis can be pulled back from 
downstairs, we can even say that $B'$ consists of one-forms  pulled back from $X$. 

To choose the modifications $L'_i$ (for $i=1,2$) locally at $P'$ we have to find a meromorphic change
of basis for the bundle $U$, which diagonalizes $\varphi ^{\ast}\theta |_U$.  The eigenforms of 
the matrix are $\pm (z')^{5}dz'$ and we can choose eigenvectors
$$
e_{\pm} := \left( 
\begin{array}{c}
z' \\
\pm (z')^{-1} \end{array} \right) .
$$
Note by calculation that 
$$
(\varphi ^{\ast}\theta |_U)  e_{\pm} =  (\pm (z')^{5}dz') e_{\pm} .
$$

Choose the line bundles $L'_1$ and $L'_2$ to be spanned by the 
meromorphic sections $e_+$ and $e_-$ of $U$. These are indeed eigen-subbundles for 
$\varphi ^{\ast}\theta$.  We just have to calculate the connection $\varphi ^{\ast} \nabla$ on
the bundle $U'=L'_1\oplus L'_2$. Which is the same as the modification of $U$ given by the meromorphic
basis $z'e_1, (z')^{-1}e_2$. 

Note first that the matrix $B'$ of one-forms pulled back from $X$
consists of one-forms which have zeros at least like $(z')^3dz'$. Thus $B'$ transported to $U'$ is still
a matrix of holomorphic one-forms so it doesn't affect our lemma. In particular we just have to consider
the transport to $U'$ of the connection $d_U$ constant with respect to the basis 
$(e_1,e_2)$ on the bundle $U$. 

Calculate
$$
d_U(a_+e_+ + a_-e_-) = 
d_U \left( 
\begin{array}{c}
(a_+ + a_-)z' \\
(a_+-a_-) (z')^{-1} \end{array} \right) 
$$
$$
=
\left( 
\begin{array}{c}
(da_+ + da_-)z' \\
(da_+-da_-) (z')^{-1} \end{array} \right) 
+
\left(
\begin{array}{c}
(a_+ + a_-)(d\log z')z' \\
-(a_+-a_-)(d\log z') (z')^{-1} \end{array} \right) 
$$
and with the notation $d_{U'}$ for the constant connection on the bundle $U'$ 
with respect to its basis $e_{\pm}$, this is equal to
$$
= d_{U'}(a_+e_+ + a_-e_-) + a_+ (d\log z') e_- + a_- (d\log z') e_+.
$$
We conclude that the connection matrix $\beta $ is, up to a holomorphic piece, just the
$2\times 2$ matrix
$$
\left(
\begin{array}{cc}
0& (z')^{-1} \\
(z')^{-1} & 0 \end{array} \right) dz'.
$$
In particular the diagonal terms of $\beta $ are holomorphic, as desired for the lemma.

These local modifications piece together to give global modifications $L'_i$ of the $L_i$.
We have to show that the $L'_i$ are of degree zero. 

In general, given a meromorphic connection on a bundle which is a direct sum of line bundles, we can
extract its ``diagonal'' part, which in terms of a local framing compatible with the direct sum
is just the connection given by the
diagonal entries of the original connection matrix. Denote this operation by $(\,\, )_{\rm diag}$.
Note that for any diagonal connection $\nabla '$ and meromorphic endomorphism-valued one-form $\beta$,
the diagonal connection is given by $(\nabla '+\beta )_{\rm diag}=\nabla ' + \beta _{\rm diag}$
where $\beta _{\rm diag}$ is the matrix of diagonal entries of $\beta$. 

Setting $E':= \bigoplus L'_i$ we have a meromorphic map $\psi : E \rightarrow E'$. We obtain a 
meromorphic connection $\psi \circ \varphi ^{\ast} \nabla \circ \psi ^{-1}$ on $E'$, and by the above
choice of $L'_i$ the associated diagonal connection is holomorphic at the singularities. On the 
other hand, $\psi \circ \varphi ^{\ast} \nabla \circ \psi ^{-1}$ is holomorphic away from the singularities,
so its diagonal part is holomorphic there too. Therefore the global diagonal connection 
$(\psi \circ \varphi ^{\ast} \nabla \circ \psi ^{-1})_{\rm diag}$ on $\bigoplus L'_i$ is holomorphic.
This proves that the $L'_i$ are of degree zero.  In particular, our choice of modification is allowable for
the argument given at the start of the  present section. This proves the lemma. 
\end{proof}

{\em Remarks:} 
\newline
(i)\,\, The above proof gives further information: the only terms with poles in 
the matrix $\beta$ are the off-diagonal terms corresponding to the two eigenvalues which
came together originally downstairs in $X$; and these terms have exactly  logarithmic (i.e. first-order)
poles with residue $1$. This information might be useful in trying to improve the current
results in order to obtain precise expansions at the singularities of the Laplace transform of the
monodromy. 
\newline
(ii)\,\, This gauge transformation is probably not new, but I don't currently have a good reference.
It looks related to \cite{Kamvissis}, \cite{Woodhouse} and \cite{SanguinettiWoodhouse}, and indeed
may go back to \cite{DDthesis} \cite{Voros}.  
\newline
(iii)\,\, The fact that we had to go to a covering whose ramification power is divisible by $4$
rather than just $2$ (as would be sufficient for diagonalizing $\theta$) is somewhat 
mysterious; it probably indicates that we (or some of us at least) don't fully understand what
is going on here. 

\medskip

Let $\beta ^{\rm diag}$ denote the matrix of diagonal entries of $\beta$.
Let $Z=\widetilde{Y}$ be the universal covering. Over $Z$ we can use the diagonal connection 
$\nabla ' + \beta ^{\rm diag}$ to trivialize
$$
E'|_Z \cong \Oo _Z^r.
$$
With respect to this trivialization , our family now has the form of a family of connections 
$$
\{ (\Oo _Z^r, d + B + tA ) \}
$$
where $A$ (corresponding to the pullback of $\theta '$ to $Z$) is the diagonal matrix whose entries
are the pullbacks of the $\alpha _i$; and where $B$ is a matrix whose diagonal entries are zero, and
whose off-diagonal entries are meromorphic with poles at the points lying over singular points for
$\theta$. 

We can now apply the method developped in \cite{abmade} to this family of connections. 
Note that it is important to know that the diagonal entries of $A$ come from forms on the compact
Riemann surface $Y$; on the other hand the fact that $B$ is only defined over the 
universal covering $Z$ is not a problem. The next two sections will constitute a brief discussion of
how the method of \cite{abmade} works; however the reader is refered back there for the full details.

\section{Laplace transform of the monodromy operators}
\label{laplace}

We now look at a family of connections of the form $d+B+tA$ on the trivial bundle $\Oo ^r$ on
the universal covering $Z$ of the
ramified cover $Y$, where $A$ is a diagonal matrix with one-forms $\alpha _i$ along the diagonal,
and $B$ is a matrix of meromorphic one-forms with zeros on the diagonal. We assume that the poles
of $B$ are at points $P\in \Rr$ coming from the original singular points of the Higgs field $\theta$
on $X$. We make no assumption about the order of poles, in spite of
the additional information given by Remark (i) after the proof of Lemma \ref{modif} above.

Assume that $p$ and $q$ are two points in $Z$, not on the singular points. 
Choose a path $\gamma$ from $p$ to $q$ not passing through the
singular points. We obtain the {\em transport matrix} $m(t)$ for continuing solutions of
the ordinary differential equation $(d+B+tA)f=0$ from $p$ to $q$ along the path 
$\gamma$.  Note that $m(t)$ is a holomorphic $r\times r$-matrix-valued function defined for all $t\in \cc$.

Denote by $Z^{\ast}$  (resp. $Z^{\epsilon}$) the complement of the inverse image of $\Rr$ 
(resp. the complement of the union of open discs of radius $\epsilon$ around points in the inverse image of $\Rr$).
The poles of $B$ force us to work in $Z^{\ast}$ rather than $Z$, and in the course of the argument an
$\epsilon$ will be chosen so that we really work in $Z^{\epsilon}$. Actually it turns out that the
fact of staying inside these regions will be guaranteed by our choice of vector fields, so we don't need
to worry about any modification of the procedure of \cite{abmade} because of this difference. 

Recall that after a gauge transformation and an expansion 
as a sum of iterated integrals, we obtain a formula for the transport matrix. One way of thinking
of this formula is to look at the transport for the connection $d+sB +tA$ and expand in a Taylor series in 
$s$ about the point $s=0$, then evaluate at $s=1$. The terms in the expansion are the higher derivatives 
in $s$, at $s=0$, which are functions of $t$. A concrete derivation of the formula is given in \cite{abmade}.
It says
$$
m(t) =\sum _I \int _{\eta _I}b_Ie^{tg_I}
$$
where:
\newline
---the sum is taken over multi-indices of the form $I= (i_0,i_1,\ldots , i_k)$ where we note $k=|I|$;
\newline
---for a multi-index $I$ we denote by $Z^{\ast}_I$ the product of $k=|I|$ factors $Z\times \ldots \times Z$;
\newline
---in $Z^{\ast}_I$ we have a cycle 
$$
\eta _I
:= \{ (\gamma (t_1),\ldots , \gamma (t_k)) \}
$$
for $0\leq t_1 \leq \ldots \leq t_k \leq 1$ where $\gamma$ is viewed as a path parametrized by $t\in [0,1]$;
\newline
---the cycle $\eta _I$ should be thought of as representing a class in a relative homology group of $Z^{\ast}_I$
relative to the simplex formed by points where $Z^{\ast}_I=z_{i+1}$ or at the ends $z_1=p$ or $z_k=q$;
\newline
---the matrix $B$ leads to a (now meromorphic) matrix-valued $k$-form $b_I$  on $Z^{\ast}_I$ defined as follows: if the 
entries
of $B$ are denoted $b_{ij}(z)dz$ then 
$$
b_I = b_{i_ki_{k-1}}(z_k)dz_k \wedge \ldots \wedge  b_{i_1i_0}(z_1)dz_1 {\bf e}_{i_ki_0}
$$
where ${\bf e}_{i_ki_0}$ denotes the elementary matrix with zeros everywhere except for a $1$ in the $i_ki_0$ place;
\newline
---and finally $g_I$ is a holomorphic function $Z^{\ast}_I \rightarrow \cc$ defined by integrating the one-forms
$\alpha _i$ as follows:
$$
g_I (z_1,\ldots , z_k) = \int _p^{z_1} \alpha _{i_0} + \ldots + \int _{z_k}^q \alpha _{i_k} .
$$

The terms in the above expression correspond to what Ecalle calls the {\em elementary monomials} ${\rm som}$,
see his article in \cite{Boutet}. 

The fact that $b_I$ is meromorphic rather than holomorphic is the only difference between our present
situation and the situation of \cite{abmade}. 
Note that because our path $\gamma$ misses the singular points and thus the poles of $B$, the cycle
$\eta_I$ is supported away from the poles of $b_I$. We will be applying essentially the same technique
of moving the cycle of integration $\eta$, but we need to do additional work to make sure it stays
away from the poles of $b_I$.  

It is useful to have the formula
$$
g_I (z_1,\ldots , z_k) = g_{i_0i_1}(z_1) + \ldots + g_{i_{k-1}i_k}(z_k)+ \int _p^q\alpha _{i_k},
$$
where
$$
g_{ij}(z):= \int _p^z \alpha _i-\alpha _j .
$$

Our formula for $m$ gives a preliminary bound of the form
$$
|m(t)| \leq Ce^a|t| .
$$
Indeed, along the path $\gamma$ the one-forms $b_{ij}$ are bounded, so 
$$
|b_I| \leq C^k
$$
on $\eta _I$; also we have a bound $| g_I(z) |\leq a$ for $z\in \eta _I$, uniform in $I$; and finally
the cycle of integration $\eta _I$ has size $(k!)^{-1}$. Putting these together gives the bound for $m(t)$
(and, incidentally, shows why the formula for $m$ converged in the first place).

Recall now that the {\em Laplace transform} of a function $m(t)$ which satisfies a bound such as the above,
is by definition the integral 
$$
f(\zeta ):= \int _0^{\infty} m(t)e^{-\zeta t}dt
$$
where $\zeta\in \cc$ with $|\zeta |> a$ and the path of integration is taken in a suitably chosen direction so that
the integrand is rapidly decreasing at infinity. In our case since $m(t)$ is a matrix, $f(\zeta )$ is also 
a matrix. We can recover $m(t)$ by the inverse transform
$$
m(t)=\frac{1}{2\pi i}\oint f(\zeta )e^{\zeta t}d\zeta 
$$
with the integral being taken over a loop going around once counterclockwise in the region $|\zeta |>a$.

The singularities of $f(\zeta )$ are directly related to the asymptotic behavior of $m(t)$. This is a classical subject
which we discuss a little bit more in \S \ref{conclusion}. One can note for example that by the inverse transform, there exist
functions $m(t)$ satisfying the preliminary bound $|m(t)| \leq Ce^a|t|$ but such that the Laplace transforms 
$f(\zeta )$ have arbitrarily bad singularities in the region $|\zeta |\leq a$. Thus getting any nontrivial restrictions
on the singularities of $f$ amounts to a restriction on which types of functions $m(t)$ can occur.

In our case, the expansion formula for $m(t)$ leads to a similar formula
for the Laplace tranform, which we state as a lemma. Define the image support of a collection $\eta = \{ \eta _I \}$
by the collection of functions $g=\{ g_I \}$ to be
the closure of the union of the images of the component pieces:
$$
g(\eta ):= \overline{\bigcup _Ig_I(| \eta _I |)}\subset \cc ,
$$
where $| \eta _I | \subset Z_I$ is the usual support of the chain $\eta _I$.

\begin{lemma}
\label{formula}
With the functions $g_I$, the forms $b_I$, and the chains $\eta _I$ intervening above, 
for any $\zeta$ in the complement of the region $g(\eta )$ the formula
$$
f(\zeta ) = \sum _I \int _{\eta _I}\frac{b_I}{g_I-\zeta} 
$$
converges, and gives an analytic continuation of the Laplace transform in the (unique) unbounded connected component
of the complement of $g(\eta )$. 
\end{lemma}
\begin{proof}
The convergence comes from the same bounds on $b_I$ and the size of $\eta _I$ which allowed us to bound $m$.
The fact that this formula gives the Laplace transform is an exercise in complex path integrals. 
\end{proof}

The terms in this expansion correspond to Ecalle's elementary monomials ``${\rm soc}$'' in \cite{Boutet}. 

A first approach would be to try to move the path $\gamma$ so as to move the union of images 
$g(\eta )$ and analytically continue $f$ to a larger region. This works quite well for rank $2$,
where one can get an analytic continuation to a large region meeting the singularities \cite{DDthesis}. 
In higher rank, the $3\times 3$ example
at the end of \cite{abmade} shows that this approach cannot be optimal. In fact, we should instead move
each cycle of integration $\eta_I$ individually. Unfortunately this has to be done with great care in order
to maintain control of the sizes of the individual terms so that the infinite sum over $I$ still converges. 

Now we get to the  main definition. It is a weak version of resurgence, see \cite{DDP}, \cite{CNP}. 

\begin{definition}
\label{lfb}
A function such as $f(\zeta )$ defined on
$|\zeta |>a$ is said to have an {\em analytic continuation with locally finite branching} if for every 
$M >0$ there is a finite set of points $S_M\subset \cc$ such that if $\sigma$ is any piecewise linear
path in $\cc -S_M$
starting at a point where $|\zeta |>a$ and such that the length of $\sigma$ is $\leq M$, then $f(\zeta )$
can be analytically continued along $\sigma$. 
\end{definition}

And the statement of the main theorem.

\begin{theorem}
\label{main}
Suppose $m(t)$ is the transport matrix from $p$ to $q$ for a family of connections on the trivial bundle 
$\Oo_Z^r$ of the form $\{ d+B + tA\}$. Suppose that $A$ is diagonal with one-forms $\alpha _i$, coming
from the pullback of a general Higgs field $\theta$ over the original curve $X$, and suppose that $B$ is
a meromorphic matrix of one-forms with poles only at points lying over the singular points of $\theta$.
Let $f(\zeta )$ denote the Laplace transform of $m(t)$. Then $f$ has an analytic continuation with locally
finite branching. 
\end{theorem}

Most of the remainder of these notes is devoted to explaining the proof.

\section{Analytic continuation of the Laplace transform}
\label{analytic}

We now recall the basic method of \cite{abmade} for moving the cycles $\eta _I$ to obtain an analytic
continuation of $f(\zeta )$. We refer there for most details and concentrate here just on stating what the end
result is. Still we need a minimal amount of notation. Before starting we should refer to \cite{DelabaereHowls}
(and the references therein) for an extensive discussion of this process for each individual integral in the sum, 
including numerical results on how the singularities of the analytic continuations determine the asymptotics of the
pre-transformed integrals.

We work with {\em pro-chains} which are formal sums
of the form $\eta = \sum _I\eta _I$ of chains on the $Z^{\ast}_I$.
We have a boundary operator denoted $\partial + A$ where
$\partial$ is the usual boundary operator on each $\eta _I$ individually,
and $A$ (different from the matrix of one-forms considered above) is a signed sum of face maps corresponding
to the inclusions $Z^{\ast}_{I'}\rightarrow Z^{\ast}_I$ obtained when some $z_i=z_{i+1}$. Our original pro-chain of integration
in the integral expansion satisfies
$(\partial + A)\eta =0$. We can write the expansion formula of Lemma \ref{formula}
as an integral over the pro-chain $\eta = \sum _I\eta _I$,
$$
f(\zeta ) = \int _{\eta} \frac{b}{g-\zeta}
$$
where $b$ is the collection of forms $b_I$ on $Z^{\ast}_I$ and $g$ is collection of functions $g_I$. 
Such a formula is of course subject to the condition that the infinite sum of integrals converges.

In a formal way (i.e. element-by-element in the infinite sums implicit in the above notation),
if we add to $\eta$ a boundary term of the form $(\partial + A)\kappa$
then the integral doesn't change:
$$
\int _{\eta + (\partial + A)\kappa } \frac{b}{g-\zeta} = \int _{\eta} \frac{b}{g-\zeta} .
$$
This again is subject to the condition that the infinite sums on both sides converge absolutely
and in fact that the individual terms in the rearrangement (i.e. separating $\partial$ and $A$) 
converge absolutely. Whenever we use this, we will be refering (perhaps without mentionning it further)
to the work on convergence which
was done in \cite{abmade}. 

Our analytic continuation procedure rests upon consideration of the locations of the images 
by the function $g$, of the pro-chains of integration. Recall the notation
$$
g(\eta ):= \overline{\bigcup _Ig_I (|\eta _I|)}
$$
where $|\eta _I|$ is the usual support of the chain $\eta _I$. 

If $f$ is defined by the right-hand integral over $\eta$ in a neighborhood of a point $\zeta_0$,
meaning that the image $g(\eta )$ misses an open neighborhood of $\zeta_0$, and if 
the image $g(\eta + (\partial + A)\kappa )$ misses an entire segment going from $\zeta_0$ to 
$\zeta _1$, then the integral over $\eta + (\partial + A)\kappa$ defines an analytic continuation
of $f$ along the segment. The procedure can be repeated with $\eta$ replaced by 
$\eta + (\partial + A)\kappa$.

At this point we let our notation slide a little bit, and denote by $\eta$ any pro-chain 
which would be obtained from the original chain of integration by a sequence of modifications of the
kind we are presently considering, such that the integral over $\eta$ serves to define an analytic continuation of 
$f(\zeta )$ to a neighborhood of a point $\zeta _0\in \cc$. The original pro-chain $\eta$ of Lemma \ref{formula}
is the initial case. Our assumption on $\eta$ says among other things that
the image $g(\eta )$ doesn't meet a disc around $\zeta _0$. Fix a line segment $S$ going from $\zeta _0$ to
another point $\zeta _1$; we would like to continue $f$ in a neighborhood of $S$. By making a rotation 
in the complex plane (which can be seen as a rotation of the original Higgs field) we may without loss
of generality assume that the segment $S$ is parallel to the real axis and the real part of $z_1$ is 
smaller than the real part of $\zeta _0$. Let $u$ be a cut-off function for a neighborhood of $S$ and
write
$$
\eta = \eta ' + \eta '' ,\;\;\; \eta ' = g^{\ast}(u) \cdot \eta .
$$
We will apply the method of \cite{abmade} to move the piece $\eta '$ (this piece corresponds to what was
called $\eta$ in Chapter 4 of \cite{abmade}). 

The first step is to choose flows.
This corresponds to Chapter 3 of \cite{abmade}. 
 In our case, we will use flows along vector fields $W_{ij}$ which
are $\mathcal{C}^{\infty}$ multiples of the gradient vector fields of the real parts $\Re g_{ij}$. 
To link up with the terminology of \cite{abmade}, these vector fields determine flowing functions
$f_{ij}(z,t)$ (for $z\in Z$ and $t\in \rr ^+$ taking values in $Z$) by the equations
$$
\frac{\partial}{\partial t}f_{ij}(z,t)= W_{ij}(f_{ij}(z,t)),\;\;\; f_{ij}(z,0) = z.
$$
Note that this choice is considerably simpler than that of \cite{abmade}.
The choice of vector fields will be discussed in detail below, and will in particular be
subject to the following constraints.

\begin{condition}
\label{constraints}
(i)\, the vector fields $W_{ij}$ are lifts to $Z$ of vector fields defined on the compact surface $Y$; 
\newline
(ii)\, the differential $d\Re g_{ij}$ applied to $W_{ij}$ at any point, is a real number $\leq 0$;  
\newline
(iii)\, there exists $\epsilon$ such that the flows preserve $Z^{\epsilon}$
i.e. the vector fields $W_{ij}$ are identically zero in the discs of radius $\epsilon$ around the
singular points; and 
\newline
(iv)\, the $W_{ii}$ are identically zero.
\end{condition}

The flows given by our vector fields
lead to a number of operators $F$, $K$ and $H$ defined as in Chapters 4 and 5 of \cite{abmade}. These
give pro-chains
$$
F\tau = \sum _{r,s} F(-KA)^rH(AK)^s\eta ',
$$
$$
F\psi = \sum _{r} F(-KA)^rK(\partial + A)\eta ',
$$
$$
FK\varphi  = \sum _{r} FK(AK)^r\eta '.
$$
The reader can get a fairly good idea of these definitions from our discussion of the points on
$|F\tau |$ in \S \ref{description} below.

\begin{lemma}
\label{assuming}
With these notations, and assuming that the vector fields satisfy the constraints marked above,
we can write
$$
\eta + (\partial + A)FK\varphi = \eta '' + F\tau - F\psi .
$$
On the right, the images $g(\eta '')$ and $g(F\psi )$ miss a neighborhood of the segment $S$. 
{\em Assuming} we can show that the image $g(F\tau )$ also misses a neighborhood of the segment $S$,
then 
$$
f(\zeta ) = \int _{\eta + (\partial + A)FK\varphi} \frac{b}{g-\zeta}
$$
gives an analytic continuation of $f$ from $\zeta _0$ to $\zeta _1$ along the segment $S$.
\end{lemma}
\begin{proof}
The operator $K$ corresponds to applying the flows defined by $W_{ij}$  in the various coordinates. This has
the effect of decreasing the real part $\Re g$. The fact that in our case we use flows along vector fields
which are positive real multiples of $-\grad \Re g_{ij}$ (this is the second of the
constraints on $W_{ij}$) implies that  the flows strictly respect the imaginary part of $g$.
This differs from the case of \cite{abmade} and means we can avoid discussion of ``angular sectors'' such as
on pages 52-53 there. Thus, in our case, when we apply a flow to a point, the new point has the same value of $\Im g$,
and the real part $\Re g$ is decreased. 

The operator $F$ is related to the use of buffers; we refer to \cite{abmade} for that discussion and heretofore ignore it.
The operator $A$ is the boundary operator discussed above; and the operator $H$ is just the result of doing
the flows $K$ after unit time. In particular, $A$ doesn't affect the value of $g$. And $H$ decreases $\Re g$
while fixing $\Im g$ just as $K$ did (this point will perhaps become clearer with the explicit description
of points in the supports of $F\tau$ and $FK\varphi$ in the next section). 

The proof of the first formula is the same as in \cite{abmade} Lemma 4.4, and we refer there for it. 

To show that the supports of $g(\eta '')$ and $g(F\psi )$ miss a neighborhood of $S$, it is useful to 
be a little bit more precise about the neighborhoods which are involved. 
Let $N_1$ be the support of $u$, which is a neigborhood of $S$ (we assume it is convex), 
and let $N_2$ be the support of $du$
which is an oval going around $S$ but not touching it. Let $N_3$ be the neighborhood of $S$ where $u$  is identically $1$.
Let $D$ be a disc around $\zeta _0$, such that
$g(\eta )$ misses $D$, and which we
may assume has radius bigger than the width of $N_1$. Then 
$$
g(\eta ') \subset N_1 - (N_1\cap D),
$$
$$
g(\eta '') \subset \cc  - (N_3\cup D),
$$
and $(\partial A)\eta ' = - (\partial A)\eta ''$ with
$$
g((\partial A)\eta ' ) \subset N_2-(N_2\cap D). 
$$
In particular the support of $g(\eta '')$ misses the neighborhood $N_3$ of $S$. 
Also, given that the boundary term $(\partial + A)\eta '$ is supported 
in the $U$-shaped region $N_2$, the effect of our operators on $\Re g$ and $\Im g$ described above
implies that $g(F\psi )$ is supported away from $N_3$. This completes the proof of the second statement of the lemma. 

For the last statement, assume that we have chosen things such that the support of $g(F\tau )$ also misses
$S$. This is certainly what we hope, because of the inclusion of the operator $H$ applying all the flows for unit time.
The only possible problem would be if we get too close to singular points; that is the technical difficulty which is
to be treated in the remainder of the paper. For now, we assume that this is done. 

Formally speaking, the first equation of the lemma means that 
$$
\int _{\eta} \frac{b}{g-\zeta} = 
\int _{\eta '' + F\tau - F\psi } \frac{b}{g-\zeta} .
$$
By our starting assumption $f(\zeta )$ is defined by the integral on
the left, in a neighborhood of $\zeta _0$. On the other hand, 
the integral on the right defines an analytic continuation along the segment $S$. 

An important part of justifying the argument of the preceding paragraph (and indeed, of showing
that the integral on the right is convergent) is to bound the sizes and numbers of all the chains
appearing here. This was done in \cite{abmade}.

The only difference in our present case is the poles in the integrand $b$.
However,  thanks to the third constraint on the vector fields $W_{ij}$,
everything takes place in $Z^{\epsilon}_I := Z^{\epsilon} \times \ldots \times Z^{\epsilon}$,
and on $Z^{\epsilon}$ there is a uniform bound on the size of $b_{ij}$. 
Also, everything takes place inside a relatively compact subset of $Z$, see \S \ref{remoteness}.
Thus the integrand in
the multivariable integral is bounded by 
$$
\sup _{Z^{\epsilon}_I}|b|\leq C^k
$$
for $k=|I|$. With this information the remainder of the argument of \cite{abmade} 
works identically the same way (it is too lengthy to recall here). 
This justifies the formal argument of two paragraphs ago and completes the proof of the lemma. 
\end{proof}

{\em Remark:} It is clear from the end of the proof that the bounds depend on $\epsilon$,
which in turn will depend on how close we want to get to a singularity. This is the root of why
we don't get any good information about the order of growth of the Laplace transform at its singularities.

\section{Description of cells using trees}
\label{description}

As was used in \cite{abmade}, the chains defined above can be expressed as sums of cells. 
We are most interested in the chain $F\tau$
although what we say also applies to the other ones such as $FK\varphi$. These chains are  unions of cells which have the
form of a family of cubes parametrized by points in one of the original cells $\eta '_I$. We call these
things just cubes. In the cubes which occur the points are parametrized by ``trees'' furnished with lots of
additional information. 
\footnote{The occurence of trees here is certainly related to and probably the same as
Ecalle's notions of {\em (co)mould (co)arborescent} cf \cite{EcalleAbordable}.
In another direction, John Conway pointed out at the time of \cite{abmade} that cubes parametrized by trees in this way
glue together into Stasheff polytopes.
I didn't know what those were at the time, but retrospectively this still remains mysterious since we are
dealing with representations of the fundamental group and it isn't clear what that 
has to do with homotopy-associativity. This is certainly a good subject for further thought.}
We make this precise as follows: a {\em furnished tree} is:
\newline
---a binary
planar tree $T$ sandwiched between a top horizontal line and a bottom horizontal line;
\newline
---with leftmost and rightmost vertical strands whose edges are called the {\em side edges};
\newline
---for each top vertex of the tree (i.e. where an edge meets the top horizontal line)
we should specify a point $z\in Z^{\ast}$ (the point corresponding to the left resp. right side edge is $p$
resp. $q$);
\newline
---for each region in the complement of the tree between the top and bottom
horizontal lines and between the side edges we should specify an index, so that each (non-side) edge of the tree is provided with left and
right indices which will be denoted $i_e$ and $j_e$ below; and
\newline
---each edge $e$ is assigned a ``length'' $s(e)\in [0,1]$.

Suppose $T$ is a furnished tree. By looking at the indices assigned to the regions meeting the top and bottom
horizontal lines we obtain multi-indices $I^{\rm top}$ and $I^{\rm bot}$, so 
the collection of points $(z_1,\ldots , z_k)$ attached to the
top vertices gives a point $z^{\rm top}\in Z^{\ast}_{I^{\rm top}}$. 

We can now explain how a furnished tree leads to a point $z^{\rm bot} \in Z_{I^{\rm bot}(T)}$.
This depends on a 
choice of vector fields $W_{ij}$ for each pair of indices $i,j$, which we now assume as having been made. 
A {\em flowing map} 
$\Phi : T \rightarrow Z$ is a map from the topological realization of the tree, into $Z$,
satisfying the following properties: 
\newline
(i)\, if $v$ is a top vertex which is assigned a point $z$ in the information
contained in $T$, then $\Phi (v)=z$;  
\newline
(ii)\, the side edges are mapped by constant maps to the points $p$ or $q$ respectively; and
\newline
(iii)\, if $e$ is an edge with left and right indices $i_e$ and $j_e$ and with initial vertex $v$ and
terminal vertex $v'$, then $\Phi (e)$ is the flow curve for flowing along the vector field $W_{i_ej_e}$
from $\Phi (v)$ to $\Phi (v')$, where the flow is done for time $s=s(e)$. This determines
$\Phi (v')$ as a function of $\Phi (v)$ and the information in the tree. Thus by recursion we determine
the $\Phi (v)$ for all vertices, as well as the paths $\Phi (e)$ for the edges $e$ (the map $\Phi$ on 
the edges is only well determined up to reparametrization
because we don't fix a parametrization of the edges; the length $s$
is abstract, since it is convenient to picture even edges assigned $s=0$ as being actual edges). 

For a given choice of vector fields $W_{ij}$ and of information attached to the tree $T$, the flowing  map
exists and is  unique. This determines a point given by the values $z$ at the bottom vertices,
$$
z^{\rm bot}(W,T)\in Z_{I^{\rm bot}}.
$$

Now go back to the situation of the previous section. Starting from a chain $\eta '$ we obtained
a chain $F\tau$. 

\begin{lemma}
\label{points}
The points in the support of $F\tau$ are described as the $z^{\rm bot}(W,T)$, where $W=\{ W_{ij}\}$ is
the collection of vector fields used to define the flows $K$ and $H$, and where $T$ is a furnished tree such that
$z^{\rm top}(T)$ is in the support of $\eta '$ and satisfying the following auxiliary condition:
\newline
(*)\, there exists 
(up to reparametrization of the planar embedding) a horizontal line
which cuts the tree along a sequence of edges, such that all of these edges are assigned the fixed length value
$s=1$.
\end{lemma}
\begin{proof}
See \cite{abmade}, pages 54-55. The auxiliary condition comes from the term $H$ in the formula for $\tau$. 
\end{proof}

{\em Remark:} For the chain $FK\varphi$ the same statement holds except that the furnished trees $T$ 
might not necessarily satisfy the auxiliary condition. 

We finish this section by pointing out the relationship between $g(z^{\rm top})$ and
$g(z^{\rm bot})$. This is the key point in our discussion, because $z^{\rm top}$ is the input point
coming from the chain $\eta '$ and $z^{\rm bot}$ is the output point which goes into the resulting
chain $F\tau$. We want to prove that the real part of
$g(z^{\rm bot})$ can be moved down past the end of the segment $S$. 

\begin{lemma}
\label{tree}
If $T$ is a furnished tree and $W$ a choice of vector fields, then 
$$
g(z^{\rm bot}(T))= g(z^{\rm top}(T)) + \sum _e \int _{\Phi (e)} dg_{i_ej_e},
$$
In particular if $W$ saitsfies Condition \ref{constraints} then 
$$
g(z^{\rm bot}(T))- g(z^{\rm top}(T)) \in \rr _{\leq 0}.
$$
\end{lemma}
\begin{proof}
If $e$ is an edge of $T$ and $s'\in [0,s(e)]$ then we can define the tree $T'$ obtained by 
{\em pruning} $T$ at $(e,s')$. This is obtained by cutting off everything below $e$ and
sending the bottom vertex of $e$ to the line at the bottom. The indices associated to regions in the
complement follow accordingly. Finally we set $s(e):=s'$ in the new tree $T'$. 

Suppose for the same edge $e$ we also pick $s'' \in [s',s(e)]$. Then we obtain a different pruning denoted
$T''$ (which has almost all the same information except for the length of the edge $e$). 
Let $v'$ (resp. $v''$) denote the bottom vertices corresponding to $e$ in the trees $T'$ (resp. $T''$).
Let $\Phi '$ (resp. $\Phi ''$) denote the flowing map for $T'$ (resp. $T''$). These coincide
and coincide with $\Phi$ on the parts of the trees that are in common (the unpruned parts).
We have
$$
g(z^{\rm bot}(T'')) = g(z^{\rm bot}(T')) + \int _{\Phi '(v')}^{\Phi ''(v'')}dg_{i_ej_e}.
$$
Note that the segment of $\Phi (e)$ going from $\Phi '(v')$ to $\Phi ''(v'')$ 
is a flow curve for the vector field $W_{i_ej_e}$, and it flows for time $s''-s'$. 

If we prune at an edge $e$ with $s'=s(e)$ then it amounts to cutting off the tree at the lower vertex of $e$.
If furthermore all of the length vectors assigned to edges below $e$ are $0$, then 
$g(z^{\rm bot}(T')) = g(z^{\rm bot}(T))$. 

By recurrence we obtain the first statement in the lemma. 

Recall that one of the constraints was the condition that 
the vector fields $W_{i_ej_e}$ be negative multiples of the gradient vector
fields for the real functions $\Re g_{i_ej_e}$.
With this condition we get that the integral of $dg_{i_ej_e}$ along a flow curve for
$W_{i_ej_e}$ is a negative real number, so this gives at each stage of the recurrence
$$
g(z^{\rm bot}(T'')) - g(z^{\rm bot}(T')) \in \rr _{\leq 0}.
$$
Putting these together gives the second statement of the lemma.
\end{proof}

There is also another way to prune a tree: if $e$ is an edge such that $i_e=j_e$ then we can cut off $e$ and
all of the edges below it, and consolidate the two edges above and to the side of $e$ into one edge.
The only difficulty here is that the consolidated edge might have total length $>1$ but this doesn't affect
the remainder of our argument (since at this point we can ignore questions about the sizes of the cells). 
Let $T'$ denote the pruned tree obtained in this way. We again have
$$
g(z^{\rm bot}(T)) - g(z^{\rm bot}(T')) \in \rr _{\leq 0}.
$$

In general we will be trying to show for the trees which arise in $F\tau$, that the real part of
$g(z^{\rm bot}(T))$ is small enough. If we can show it for $T'$ then it follows also for $T$. In this way
we can reduce for the remainder of the argument, to the case where $i_e\neq j_e$ for all edges of $T$.
This is the content of the following lemma. For its statement, recall the neighborhood 
$S\subset N_1$ appearing in the proof of Lemma \ref{assuming}.

\begin{lemma}
\label{redux}
Let $|S|=\zeta _0-\zeta _1$ denote the length of the segment along which we want to continue $f$. 
In order to show that the image $g(F\tau )$ misses a neighborhood, say $N_1$, of the segment $S$
it suffices
to choose our vector fields $W$ (satisfying Condition \ref{constraints})
so that if $T$ is any furnished tree satisfying:
\newline
(i)\, the auxiliary condition
(*) of Lemma \ref{points};
\newline
(ii)\, that $i_e\neq j_e$ for all edges $e$ of $T$; and
\newline
(iii)\, that $z^{\rm top}(T)$ is in the support of $\eta '$;
\newline 
then 
$g(z^{\rm bot}(T))$ lies outside of our neighborhood $N_1$ of $S$.
\end{lemma}
\begin{proof}
Assume that we have chosen the vector fields to give the reduced condition of this statement.
Suppose $z$ is a point on the support of $F\tau$. Then there is a furnished tree $T^1$ as in Lemma \ref{points}
such that $z=z^{\rm bot}(T^1)$ and such that $z^{\rm top}(T^1)$ is on the support of $\eta '$.
Let $T:= (T^1)'$ be the pruning of $T^1$ described directly above. It still satisfies (i), i.e. the
condition (*) of Lemma \ref{points}, and by the pruning process it automatically satisfies (ii). Also 
$z^{\rm top}(T) =z^{\rm top}(T^1)$ is on the support of $\eta '$, so our condition gives that
$g(z^{\rm bot}(T))$ lies outside of $N_1$. On the other hand, 
$$
g(z^{\rm bot}(T^1)) - g(z^{\rm bot}((T^1)')) \in \rr _{\leq 0}, \;\;\; 
g(z^{\rm bot}(T))- g(z^{\rm top}(T)) \in \rr _{\leq 0}.
$$
Thus $g(z)= g(z^{\rm bot}(T^1))$, then $g(z^{\rm bot}(T)) =g(z^{\rm bot}((T^1)'))$,
and then $g(z^{\rm top}(T))$ lie in order on a line segment parallel to the real axis. 
Given that $g(z^{\rm top}(T))\in N_1$ but $g(z^{\rm bot}(T)) \not\in N_1$, and that $N_1$ is 
a convex, we obtain $g(z)\not\in N_1$ as desired. 
\end{proof}

{\em Remark:} The condition of the lemma will not be possible, of course, 
when the segment $S$ passes through a turning point. 
Finding out the conditions on $S$ to make it possible will tell us where the turning points are.

\section{Remoteness of points}
\label{remoteness}

One of the important facets of the statements of theorems \ref{main} and \ref{mainplus} is the local
finiteness of the set of singularities. We describe here briefly how this works. It reproduces
the discussion of \cite{abmade}, but with considerable simplification due to Condition \ref{constraints} (iv)
which says that when $i=j$ the flow $f_{ij}(z,t)$ is constant. 

It should be noted that the local finiteness notion \ref{lfb} is fairly strong
in that one can wind arbitrarily many times around a given singularity
for an arbitrarily small cost in terms of length of the path. In our mechanism, this is achieved by 
analytically continuing along a large number of very small segments. 

We can choose a metric $d\sigma $ on $Z^{\ast}$ (and which is a singular but finite metric
on $Z$) with the property that for any distinct
pair of indices $i\neq j$, if $\xi : [0,1]\rightarrow Z$
is a path whose derivative is a negative real multiple of $\grad \Re g_{ij}$ then
$$
\int _{\xi}d\sigma \leq \Re (g_{ij}(\xi (0)) - g_{ij}(\xi (1))).
$$

Now suppose $z = (z_i)\in Z_I$, and suppose $T$ is a binary planar tree embedded in $Z$,
with one top vertex at $p$ and whose bottom vertices are the $z_i$. Let
$$
{\bf r}_T(z):= \int _T d\sigma
$$
be the total length of the tree with respect to our metric. Define the {\em remoteness} ${\bf r}(z)$ to
be the infimum of ${\bf r}_T(z)$ over all such trees. 

\begin{lemma}
\label{remote1}
Suppose $T$ is a furnished tree, and use flows defined by vector fields satisfying Condition 
\ref{constraints} to define $z^{\rm bot}(T)$. Then 
$$
{\bf r}(z^{\rm bot}(T)) \leq {\bf r}(z^{\rm top}(T)) + g(z^{\rm top}(T)) -g(z^{\rm bot}(T)).
$$
\end{lemma}
\begin{proof}
If $T^1$ is any tree as in the definition of remoteness for $z^{\rm top}(T)$ then we can add $T$ to $T^1$
(the top vertices of $T$ being the same as the bottom vertices of $T^1$)
to obtain a tree $T^2$ as in the definition of remoteness for $z^{\rm bot}(T)$.  The formula
$$
{\bf r}_{T^2}(z^{\rm bot}(T)) \leq {\bf r}_{T^1}(z^{\rm top}(T)) + g(z^{\rm top}(T)) -g(z^{\rm bot}(T))
$$
is immediate from Lemma \ref{tree} and the property of $d\sigma$; 
use Condition \ref{constraints} (iv) to deal with edges of $T$ 
having $i_e=j_e$. 
\end{proof}

\begin{lemma}
\label{remote2}
Let $\gamma$ be a path from $p$ to $q$, which leads to the original pro-chain $\eta$ appearing in
Lemma \ref{formula}. Suppose $M_0$ is the length of $\gamma$ in the metric $d\sigma$. Then for any point
$z$  on the support of $\eta$ we have ${\bf r}(z)\leq M_0$.
\end{lemma}
\begin{proof}
For any point $z$ on the support of $\eta$, we have $z_i=\gamma (t_i)$ for $t_1\leq \ldots \leq t_k\leq 1$.
The path $\gamma$ can be considered as a tree (of total length $M_0$) 
starting at $p$ with one spine and $k$ edges of length $0$ coming off at the points $z_i$. 
\end{proof}

In our procedure for analytic continuation along a path of length $\leq M$, we obtain chains
whose support consists only of points with ${\bf r}(z)\leq M_0+2M$ (see \S \ref{proof} below). 
In particular each $z_i$ is at 
distance $\leq M_0+2M$ from $p$ with respect to $d\sigma$. Thus everything we do takes place in
a relatively compact subset of $Z$ (and concerns only a finite number of singular points $P\in Z$).

\section{Calculations of gradient flows}
\label{calculations}

We express the gradient of the real part of a holomorphic function,
as a vector field in a usual coordinate and in logarithmic coordinates.
This is of course elementary but we do the calculation just to get the
formula right. 
Suppose $z$ is a coordinate in a coordinate patch on $X$. 
The metric on $X$ may be expressed by the real-valued positive 
function
$$
h(z):= \frac{|dz |^2}{2}.
$$
Write $z=x+iy$. Note that $dx$ and $dy$ are perpendicular and have
the same length, so 
$$
h (z) = |dx|^2.
$$
The real tangent space has orthogonal
basis 
$$
\{ \frac{\partial}{\partial x}, \frac{\partial}{\partial y} \}
$$
and the formula
$$
1 = | \frac{\partial}{\partial x} \cdot dx| = h^{1/2}
| \frac{\partial}{\partial x} |
$$
yields 
$$
| \frac{\partial}{\partial x} | = h^{-1/2}.
$$
In particular an orthonormal basis for the real tangent space is
given by 
$$
\{ h^{1/2}\frac{\partial}{\partial x},
h^{1/2} \frac{\partial}{\partial y} \} .
$$
Thus we have the formula, for any function $a$:
$$
\grad a = h\frac{\partial a}{\partial x} \frac{\partial}{\partial x}
+
h\frac{\partial a}{\partial y} \frac{\partial}{\partial y} .
$$

Now suppose $g=a+ib$ is a holomorphic function (with $a,b$ real), and pose 
$f(z):= \frac{\partial g}{\partial z}$ so that $dg = f(z)dz$.
Write $f(z)=u+iv$ with $u,v$ real, and expand:
$$
(u+iv)(dx + idy) = \frac{\partial a}{\partial x} dx +
\frac{\partial a}{\partial y}dy + 
i  \frac{\partial b}{\partial x} dx +
i\frac{\partial b}{\partial y}dy.
$$
Comparing both sides we get
$$
u=\frac{\partial a}{\partial x},\;\;\;\;
v=-\frac{\partial a}{\partial y}.
$$
Note that $a=\Re g$ is the real part of $g$, so finally we have the formula
$$
\grad \Re g = h(z) \left(
(\Re \frac{\partial g}{\partial z} ) \frac{\partial}{\partial x}
- (\Im \frac{\partial g}{\partial z}) \frac{\partial}{\partial y} \right) .
$$

Suppose now that $w$ is a local coordinate at a point $P$, and consider
$$
g = a_mw^m.
$$
Let $z=-i\log w$ so $w = e^{iz}$, and writing $z=x+iy$ we have
$w= e^{ix-y}$. Then
$$
g(z) = a_m e^{imz}; \;\;\; \frac{\partial g}{\partial z} =
mia_m e^{imz}.
$$
If we write $mia_m = e^{r+is}$ then
$$
\frac{\partial g}{\partial z} = e^{r-my + i(s+mx)},
$$
so
$$
\grad \Re g =
h(z)e^{r-my}\left(
\cos (s+mx) \frac{\partial}{\partial x} , 
\sin (s+mx) \frac{\partial}{\partial y} \right) .
$$

The {\em asymptotes} are the values $x=B$ where $\cos (s+mx)=0$. At these
points, the gradient flow vector field is vertical (going either up or down,
depending on the sign of $\sin (s+mx)$). If the flow goes up, then it
stays on the vertical line until $y=\infty$. 

Note that the gradient of $\Re g$ is perpendicular to the level curves
of $\Re g$, so it is parallel to the level curves of $\Im g$.
Which is to say that the level curves of $\Im g$ are the flow lines.
This gives an idea of the dynamics of the flow. 
We have
$$
\Im g =\Im ( im^{-1}e^{(r-my)+i(s+mx)}) = -m^{-1}e^{r-my}\sin (s+mx).
$$
Thus a curve $\Im g = C$ is given by
$$
e^{-my} = \frac {-mC}{e^r\sin (s+mx)}
$$
or (noting that the sign of $C$ must be chosen so that the right hand
side is positive)
$$
y = m^{-1}r\log | \sin (s+mx) |- m^{-1}\log |mC| .
$$
In particular the level curves are all vertical translates of the same
curve; this curve $y = m^{-1}r\log |\sin (s+mx)|$ has vertical
asymptotes at the points where $\sin (s+mx)=0$. Note however that
at the asymptotes, we get $y\rightarrow -\infty$; whereas our coordinate
patch corresponds to a region $y>y_0$. Thus, every gradient flow
except for the  inbound (i.e. upward) flows directly on the asymptotes,
eventually turns around and exits the coordinate patch. This of
course corresponds to what the classical picture looks like in terms
of the original coordinate $w$. 

Also we can calculate the second derivative (which depends only on $x$
and not on which level curve we are on, since 
they are all vertical translates). Consider for example points where $\sin (s+mx)>0$.
There
$$
\frac{dy}{dx} = \frac{r\cos (s+mx)}{\sin (s+mx)}
$$
and
$$
\frac{d^2y}{dx^2} = \frac{-rm}{\sin ^2(s+mx)}
$$
In particular note that we have a uniform bound everywhere:
$$
\frac{d^2y}{dx^2} \leq -\gamma ,
$$
here with $\gamma = rm$.

Suppose  now more generally that $g$ is a holomorphic function with
Taylor expansion
$$
g = a_mw^m + a_{m+1}w^{m+1} +\ldots .
$$
Then we will get 
$$
h(z)^{-1}e^{my-r}\grad \Re g =
\left(
\cos (s+mx) \frac{\partial}{\partial x} , 
\sin (s+mx) \frac{\partial}{\partial y} \right)  + O (e^-y).
$$
In particular, the direction of the gradient flow for $g$ is determined,
up to an error term in $O(e^{-y})$, by the vector 
$(\cos (s+mx),  \sin (s+mx) )$.

The asymptotes are no longer vertical curves, but they remain in bands $x\in B_{ij,a}$.
Also we can choose $A$ in the definition of steepness, so that at non-steep parts of
the level curves we still have a bound 
$$
\frac{d^2y}{dx^2} \leq -\gamma .
$$

\section{Choice of the vector fields $W_{ij}$}
\label{choice}

The only thing left to be determined in order to fix our procedure for moving the
cycle of integration is to choose the vector fields. 
Before going further, fix a smooth metric $h$ on $Z$, for example coming from the pullback of
a smooth metric on $Y$. Use this to calculate gradients. Suppose $\epsilon$ is given. Let 
$\rho$ denote a cutoff function which is identically $0$ in the discs $D_{\epsilon /2}(P)$ (for all
points $P$ in the inverse image of $\Rr$), and is identically $1$ outside the (closed) discs $D_{\epsilon }(P)$.
Of course $\epsilon$ will be small enough that the discs don't intersect. Consider also a positive real
constant $\mu \in \rr _{>0}$. Then we put
$$
W_{ij}:= \mu \grad \Re g_{ij},
$$
and 
$$
W'_{ij}:=\rho W_{ij}.
$$
The vector fields $W'_{ij}$ satisfy Condition \ref{constraints} (with $\epsilon /2$ in place of 
$\epsilon$). We will use these vector fields for our choice of flows, and apply the criterion of
Lemma \ref{redux}. 

The point we want to make in the present section is that the flow curves for the cut-off gradient
vector field $W'_{ij}$ 
are the same as those of the true gradient flow along $W_{ij}$, up until any point where they enter
some $D_{\epsilon}(P)$. This will allow the notational simplification of looking 
at $W_{ij}$ rather than $W'_{ij}$ in the next section.

Let $\nu >0$ be the radius used to define the oval neighborhood $N_1$, i.e.\ choose
$N_1$ equal to the
set of points of distance $<\nu$ from $S$.  Once $\epsilon$ is given, choose $\mu$ large enough so
that the following property holds:

\begin{condition}
\label{mucond}
If $z(t)=f_{ij}(z_0,t)$ is a flow curve for $W_{ij}$ (for distinct indices $i\neq j$)
which never enters into any $D_{\epsilon}(P)$
flowing for $t\in [0,s]$ with $s\geq 1$, then
$$
g_{ij} (f_{ij}(z_0,s)) - g_{ij}(z_0) < \zeta _1-\zeta _0 -2\nu .
$$
Recall that $\zeta_0,\zeta _1$ were the endpoints of the segment $S$ with
$\zeta _1-\zeta _0$ a negative real number. 
\end{condition}

It is possible to choose $\mu$ (we only need to do it over a relatively compact subset of $z_0\in Z$
by the remark at the end of \S \ref{remoteness}, but in any case everything involved
is pulled back from the compact $Y$ so the choice of $\mu$ is uniform in $z_0$). 

The next lemma formalizes the following reduction: the trees which show up in Lemma \ref{redux}
have a horizontal line of edges assigned length $1$. If the flow for at least one of these edges
stays outside of all the $D_{\epsilon}(P)$ then by Condition \ref{mucond} the value of $g$ is 
decreased sufficiently to get us out of $N_1$. Thus the only case which poses a problem is when
every downward branch of the tree ends up flowing into some $D_{\epsilon}(P)$. In this case we
prune the tree at the points where it enters these  discs.

\begin{lemma}
\label{onlyprob}
Suppose $\epsilon$ is given, and $\mu$ chosen to satisfy Condition \ref{mucond}. Use the
vector fields $W'_{ij}$ to define the flows. 
In order to show that the image $g(F\tau )$ misses our neighborhood $N_1$
of the segment $S$, it suffices
to show that if $T$ is any furnished tree satisfying the following conditions:
\newline
(i)\, that $z^{\rm top}(T)$ lies on the support of $\eta '$;
\newline
(ii)\, that $i_e\neq j_e$ for any edge of $T$; 
\newline
(iii)\, that for each bottom vertex $v$ of $T$ there is a singular point $P(v)$ such that
$\Phi (v) \in D_{\epsilon} (P (v))$; and
\newline
(iv) \, that all other points of $\Phi (T)$ are outside the discs $D_{\epsilon} (P)$,
\newline
then $g(z^{\rm bot}(T))$ is not in the neighborhood $N_1$ of $S$.
\end{lemma}
\begin{proof}
Suppose $T$ is a furnished tree as in the reduction of Lemma \ref{redux}.
Prune $T$ at any point where the flowing map $\Phi$ enters into one of the closed discs $D_{\epsilon}(P)$. 
If this prunes all branches of the tree, then by an argument using \ref{tree} similar to the previous reductions,
that puts us in the case described here so we are done.
 
Thus we may assume that there is at least one branch which is not pruned. By condition \ref{redux} (i)
which is the same as Condition (*) of Lemma \ref{points}, the branch going to the bottom
has at least one edge assigned length $1$. This edge has $i_e\neq j_e$. By Condition \ref{mucond}
we have for this edge
$$
\int _{\Phi (e)}dg_{i_ej_e} < \zeta _1 -\zeta _0 -2\nu .
$$
Therefore, by the formula of Lemma \ref{tree} we have 
$$
g(z^{\rm bot}(T)) -g(z^{\rm top}(T)) < \zeta _1 -\zeta _0 -2\nu .
$$
Given that $g(z^{\rm top}(T)) \in g(\eta ') \subset N_1$ but $N_1$ is an oval with largest diameter
$2\nu + \zeta _0-\zeta _ 1$, we get $g(z^{\rm bot}(T)) \not\in N_1$. 
\end{proof}

\begin{corollary}
\label{onlyprob2}
Define the chain $F\tau$ using the vector fields $W'_{ij}$. Then, in order to show that 
$g(F\tau )$ misses $N_1$ it suffices to show that for any furnished tree $T$ satisfying 
the conditions (i)-(iv) of \ref{onlyprob} with respect to the flowing map $\Phi$ defined by the vector fields 
$W_{ij}$ (rather than $W'_{ij}$), we have 
$g(z^{\rm bot}(T))\not \in N_1$.
\end{corollary}
\begin{proof}
The two flowing maps coincide, in view of condition (iv).
\end{proof}

In view of this corollary, we can in the next section ignore the cutoff functions $\rho$ and look directly at the
gradient flows $W_{ij}$.

\section{Results on the dynamics of our flowing maps}
\label{results}

We will consider a system
of discs centered at our singular points $P$:
$$
D_{\epsilon}(P) \subset D_{\xi}(P) \subset D_u(P) \subset D_w(P).
$$
We will first fix $u$ and $w$ so that certain things are true in
a coordinate system for $D_w(P)$ (and say $u=w/2$). Then once
$u$ and $w$ are fixed we will let $\epsilon\rightarrow 0$.
Finally $\xi >\epsilon$ will be a function of $\epsilon$ with
$\xi \rightarrow 0$ when $\epsilon \rightarrow 0$.

The innermost discs $D_{\epsilon}(P)$ are those which will enter into the reduction
of Lemma \ref{onlyprob}. 
Recall that $\mu$ is chosen after $\epsilon$. 
In view of the Corollary \ref{onlyprob2}, we henceforth look
directly at the gradient flows $W_{ij}=\mu \grad \Re g_{ij}$. 

Our first lemma bounds the number of outgoing subtrees.

\begin{lemma}
\label{outgoing}
If $T$ is a furnished tree with one top edge $e$, and if $\Phi : T \rightarrow
X$ is a flowing map such that the images of all bottom vertices
are contained in some $D_{\epsilon}(P_i)$, and if
$\Phi (e)$ exits from $D_u(P)$ then $T$ contains a strand $\sigma$
such that $\Phi (\sigma )$ exits from $D_w(P)$ also. 
\end{lemma}

Our next lemma gives a normal form for any subtree which stays
entirely within $D_u(P)$. 

\begin{lemma}
\label{special}
If $T$ is a furnished tree with one edge $e$ at the top, and if $\Phi$ 
is a flowing map from $T$ into $D_u(P)\subset X$ such that
all of the bottom vertices are mapped into $D_{\epsilon}(P)$,
then the curve $\Phi (e)$ passes into $D_{\xi}(P)$, 
and flows along a vector field $W_{i_ej_e}$ in an ingoing sector 
near an ingoing curve $G_{i_ej_e}$. 
\end{lemma}

The last of our preliminary lemmas bounds the number of subtrees 
having the previous normal form. 

\begin{lemma}
\label{innerbound}
There is a number $K$ (depending on $u,w,A$ but independent of $\epsilon$, $\xi$ and $\mu$)
such that if $T$ is a furnished tree consisting of
one edge strand $\kappa$ plus a number of sub-trees coming out of
$\kappa$, and if $\Phi$ is a flowing map from 
$T$ into $D_u(P)$ with the property that all the sub-trees
coming  out of $\kappa$ are covered by Lemma \ref{special},
then there are $\leq K$ of these sub-trees.
\end{lemma}

For the proofs of these lemmas, we will use a logarithmic
coordinate system for $D_w(P)$. If $z_D$ denotes the coordinate in
the disc then we introduce $z_L = -i\log z_D$ and write $z_L=x+iy$
$z_D = e^{ix-y}$. 

The disc $D_w(P)$ is given by $y>y_0$. 

The vector fields $W_{ij}=\mu \grad \Re g_{ij}$ are approximately
equal (up to a term smaller by a factor of $O(\mu e^{-y})$) to the standard
vector fields $W'_{ij}=\mu \grad \Re  g'_{ij}$ where $g'_{ij}$ is the
leading term in the Taylor expansion for $g_{ij}$ at $P$.

Because of this, we obtain the following facts. The asymptotic directions
(which are close to vertical lines) occur in bands of the form 
$x\in B_{ij,a}$ where $B_{ij,a}\subset \rr$ are intervals 
 which can be made as small as we like
by modifying 
$y_0$.  These intervals are disjoint, except for the asymptotes of the pairs
$\{W_{ik},W_{jk}\}$ or $\{ W_{ki}, W_{kj}\}$, where $i,j$ are the
two indices attached to $P$, and $k$ is any index different from these two.
In those cases the pairs share the same values $B_{ij,a}$ and the same bands.
We say that a vector field $W_{ij}$ is {\em attached} to an interval
$B$ if $B = B_{ij,a}$. The only intervals with more than one vector field
attached to them are those described above. 

It is worth mentionning why we have this disjointness property. It is because of the non-parallel
condition on the eigenforms of $\theta$ at the singular points. The non-parallel condition
implies that the bands, which are the solutions of $s+mx =0$ modulo $\pi$, are distinct, because
the values of $s$ (which are the angular coordinates of the constants attached to the
leading terms of $g_{ij}$ as explained in the preceding section)
are different exactly because of it. Notice that the exponents $m$ are the same for all of the values
$ij$ except the two attached to the singular point; for those which are attached the value $m'$
is bigger. The non-parallel condition gives disjointness for all of the bands except the 
ones corresponding to the attached indices $ij$ and $ji$. For those, note that if we make
a general rotation of everything, the asymptotic solutions of $s+m'x$ move differently than the
solutions of $s+mx$, so those bands are disjoint from all the other ones. The general rotation
of everything corresponds to a condition that the line segments in the complex plane along which
we analytically continue, might be constrained to not be parallel to a certain finite number of directions.
This doesn't hurt our ability to analytically-continue the function. 

We can fix a number $A>0$ with the following properties: outside
of an asymptotic band for $W_{ij}$ or $W_{ji}$, the slope of the
vector $W_{ij}$ satisfies
$$
\left|
\frac{dy}{dx} (W_{ij})
\right| \leq A.
$$
Inside an asymptotic band $B$, only the vector fields $W_{ij}$ 
which are attached to $B$ can have slope bigger than $A$ or less  than $-A$.

Suppose now that $(x(t), y(t))$ is a flow along one of the vector
fields $W_{ij}$. We say that the path is {\em steep} if 
$$
\left|
\frac{dy}{dx}
\right| > A,
$$
and we say that it is {\em not steep} otherwise. We say that the path is
{\em ingoing} if $\frac{dy}{dx} > 0$ and outgoing otherwise. 
Note that with our logarithmic coordinate system, outgoing is downward
and ingoing is upward. The coordinate patch (i.e. choice of $y_0$)
and the choice of $A$ can be made so that all of the paths satisfy
the following property: 
\newline
---once the path is steep and outgoing, it remains steep and outgoing for the
remainder of the time of definition, and ends up leaving the region $y>y_0$.
\newline
This is true even though the vector field is not exactly equal to the
standard model but only close to it. 

On the other hand, the direction, i.e. the 
sign of
$\frac{dx}{dt}$ remains the same throughout the interval where the 
path is not steep. 
Call this
sign $(-1)^m$. In particular
we can think of the path as being parametrized by $x$. Define the {\em slope}
to be the signed derivative $(-1)^m\frac{dy}{dx}$. 

We have a bound, in the region where the path is not steep:
$$
\frac{d^2y}{dx^2} \leq -\gamma
$$
with $\gamma > 0$ a positive constant. Note that the second derivative is also the
variation of the slope with respect to $x$ when we go in the direction of the path. 

In particular, once the
path is outgoing it remains outgoing for the remainder of its period
of definition. This is because of the second derivative when it
is not steep, and the fact that when it becomes steep and outgoing then
it stays that way. 

We now note the {\em additive relation} for the 
vector fields at vertices of a tree. 

\begin{lemma}
\label{additive}
Suppose we are in the
situation of a flowing map $\Phi : T\rightarrow X$ defined by vector fields
$W_{ij}=\mu \grad \Re g_{ij}$. 
At any vertex $v$ of
$T$ with edges noted $e_1,e_2,e_3$ (say $e_1$ ingoing and
$e_2,e_3$ outgoing), we have three indices $i,j,k$ such that
$$
i_{e_1}=i_{e_2}=i; \; j_{e_2}=i_{e_3}=j; \; j_{e_1}=j_{e_3}=k.
$$
For the three vector fields 
$W_{ik}, W_{ij}, W_{jk}$ corresponding to the edges $e_1,e_2,e_3$ we have the relation
$$
W_{ik}(\Phi (v))=  W_{ij}(\Phi (v)) +W_{jk}(\Phi (v)).
$$
\end{lemma}
\begin{proof}
The vector fields $W_{ij}$ are all the same multiple of the
gradients $\grad \Re g_{ij}$. The fact that 
$dg_{ij} =\alpha _i -\alpha _j$ implies that 
$dg_{ij}+dg_{jk} =dg_{ik}$ giving the relation in question.
\end{proof}

\begin{proof}[Proof of Lemma \ref{outgoing}]
The disc $D_u(P)$ will be determined by $y>y_1$ for 
some $y_1$ fixed as a function
of $y_0$ (and in fact one could take $y_1=y_0+1$ for example). 
A consequence of the additive relation is that if 
$W_{ik}(\Phi (v))$ is outgoing (i.e. $\frac{dy}{dt} \leq 0$ along
this vector) then one of the other two 
$W_{ij}(\Phi (v))$ or $W_{jk}(\Phi (v))$ will also be outgoing. 
As we have noted above, if the flow along any edge is outgoing at
some point then it is outgoing for all further points. In particular
if at any point in the tree the flow  is outgoing then we can choose
a strand going down to the bottom, along which the flow is always
outgoing.  If there is an edge which crosses
out of $D_u(P)$, at the crossing point it has $\frac{dy}{dt} \leq 0$,
so we get a strand which maintains $\frac{dy}{dt} \leq 0$ as long
as it stays inside $D_w(P)$. In particular the strand cannot go back
to $D_{\epsilon}(P)$ so it must exit from $D_w(P)$ (here using the hypothesis
that any strand must end in some $D_{\epsilon}(P_i)$). This completes
the proof of Lemma \ref{outgoing}.
\end{proof}

Now we come to the proofs of Lemmas \ref{special} and \ref{innerbound}.
Fix notations $L:= -\log \epsilon$ and $L_1 := -\log \xi$. Thus we will 
let $L\rightarrow \infty$ and we have to specify $L_1$ as a function
of $L$ such that $L_1\rightarrow \infty$ too. Our discs $D_{\epsilon}(P)$
and $D_{\xi}(P)$ respectively become the regions $y>L$ and $y>L_1$.
We will specify $L_1$ as a function of $L$ so as to make the proofs
of Lemmas \ref{special} and \ref{innerbound} work.

In both lemmas, we lift the maps $\Phi$ into  maps into the coordinate
chart for the logarithmic coordinates.

\begin{proof}[Proof of Lemma \ref{special}] At any point where the flow is
not steep, the second derivative is bounded above by $-\gamma$.
In particular the flow becomes outgoing before it becomes steep again.
Furthermore, if $v$ is a vertex with indices $i,j,k$ as above,
such that the vector field 
$W_{ik}(\Phi (v))$ is not steep but is ingoing, then the additive relation
insures that one of the other two flows 
$W_{ij}(\Phi (v))$ or $W_{jk}(\Phi (v))$ has slope less than or equal to
the slope of $W_{ik}(\Phi (v))$. For this, draw a line
through the first vector, and note that one of the two other vectors has
to lie below or on the line. Note that this gives two cases:
either the new vector changes direction (i.e. the sign $(-1)^m$ changes)
and the new vector is in fact outgoing; or else the direction stays the
same and the slope decreases. Thus if $t_0$ is any point in $T$
where the flow is ingoing but not steep, then we can choose a strand $\sigma$
below $t_0$ with the property that at the end of the strand the flow becomes
outgoing; and along the strand the direction stays the same and the
second derivative satisfies 
$$
\frac{d}{dx} ((-1)^m\frac{dy}{dx}) \leq - \gamma 
$$
in a distributional sense. Then (noting by $x(t),y(t))$ the coordinates
of the image point $\Phi (t)$ for $t\in \sigma$) we have
$$
y(t) \leq y(t_0) +A (-1)^m(x(t)-x(t_0)) - \frac{\gamma}{2} (x(t)-x(t_0))^2
$$
for any $t\geq t_0$. In particular there is a number $N$ such that
$$
y(t) \leq y(t_0)+N
$$
further along the strand.
We will choose $L_1 = L - N$. 

Recall now that in the hypotheses of the lemma, we suppose that
all strands in the tree remain inside $D_u(P)$ and also finish in $D_{\epsilon}(P)$. However, we construct above
a strand which eventually becomes outgoing; therefore the strand must enter the
region corresponding to $D_{\epsilon}(P)$ before it becomes outgoing
(and notice also that it could simply stop inside this region before 
becoming outgoing, a case not mentionned above). In particular,
if there is any point $t_0$ corresponding to a non-steep ingoing flow,
or of course to any sort of outgoing flow, then we have to have
$y(t_0)+N > L$ or $y(t_0)>L_1$. 

Now we can complete the proof of the lemma. If $v$ is any vertex, such that
the incoming edge is steep and ingoing, then one of the two outgoing edges 
has to be either non-steep and ingoing, or outgoing. This is verified from
the fact that at most two different vector fields can be attached as ingoing
asymptotic vector fields for the same band $B$. From what was said above,
the bottom vertex of the first edge $e$ of the tree must satisfy 
$y(\Phi (v)) > L_1$, in other words the first edge continues all the
way until $D_{\xi}(P)$. Also the part of the edge $e$ which is
outside of $D_{\xi}(P)$ must be contained in an ingoing asymptotic band
for its vector field $W_{ij}$ and the flow is steep at all points 
of $\Phi (e)$ which are outside of $D_{\xi}(P)$. This completes the proof 
of Lemma \ref{special}.
\end{proof}

\begin{proof}[Proof of Lemma \ref{innerbound}]
Consider a vertex $v$ along $\kappa$ where a subtree in the normal 
form of Lemma \ref{special} comes off. 
Use the same notation as previously for
the edjes and indices adjoining $v$. For the sake of simplicity we
assume that $\kappa$ corresponds to the two leftmost edges $e_1$ and $e_2$
at $v$. The upper edge of the subtree is thus $e_3$ with indices $jk$.

Note from the proof of 
\ref{special} that $W_{jk}$ is ingoing and steep at $\Phi (v)$. 

As a first case, note that if the anterior edge $e_1$ of
$\kappa$ has $W_{ik}$ which is outgoing and steep, then 
the subsequent edge $e_2$ of $\kappa$ is also outgoing and steep. 
In particular at any point where $\kappa$ becomes outgoing and steep,
it remains that way and in fact will leave the region $y>y_1$ before
it goes  into any other band $B$. By looking at the possible combinatorics of
the indices one sees, even in the case of two vector fields sharing the same
band, that there can be no further normal-form vertices on $\kappa$.

In view of the previous paragraph we may restrict our attention to the places
where $\kappa$ is either not steep, or else steep but ingoing. However,
if it is steep but ingoing then again at most one vertex with a normal-form
subtree can correspond to the current band; thus at some point $\kappa$
leaves this band and must become non-steep.
On the other hand, once $\kappa$ is non-steep, 
it doesn't change to become steep and ingoing. It doesn't do this in the middle
of an edge, because of the second derivative condition.
It doesn't do it at a vertex because the edge $e_3$ which comes off is 
steep and ingoing, and a $W_{ik}$ which is not steep 
couldn't be the sum of two steep
and ingoing vectors. 

The two previous paragraphs show that we may (at the price of at most
two extra normal-form subtrees) restrict our attention to the region where
$\kappa$ is non-steep. Now one sees again from the additive relation that
if $W_{ik}$ and $W_{ij}$ are non-steep, whereas $W_{jk}$ is steep and
ingoing, then the directions of $W_{ik}$ and $W_{ij}$ must be the same. 
Indeed, if not then we would have
$W_{jk}=W_{ik}+(-W_{ij})$
which would be a sum of two vectors in the same non-steep quadrant,
so $W_{jk}$ in a steep quadrant would be impossible.

Since the sign $(-1)^m$ of $\frac{dx}{dt}$ doesn't change,
we can use $x$ to parametrize $\kappa$.
Furthermore the slope $(-1)^m\frac{dy}{dx}$ is decreasing
along $\kappa$ (note that at any vertices where a subtree in normal
form comes off, the remaining outgoing edge of $\kappa$ has a smaller
slope than the ingoing edge, because of the additive relation). 
  
In other words, 
the second derivative is distributionally less than the constant
$-\gamma$, so at some time $t$ with $|x(t)-x_0|\leq 2A /\gamma$
we get to $(-1)^m\frac{dy}{dx}\leq -A$, i.e. $\kappa$ becomes steep and
outgoing. 
We get that the non-steep part of 
the path $\kappa$ is parametrized by an interval in the
$x$-coordinate, of length $\leq 2A/\gamma $. 
There is a bound $K$ so that
 such an interval can cross (or go near) 
at most $K-2$ asymptotic bands. A band is attached to at most two pairs
of indices, but only one of these can lead correspond to a normal-form
subtree. 
Thus (counting the two we may have missed above)
the number of normal-form subtrees attached
to $\gamma$ is $\leq K$. This completes the proof of Lemma \ref{innerbound}.
\end{proof}

We  now come to the main result of this section. 
Fix $u,w$ as above, and let $L_1:= L+N$ be the function determined by
the above proofs. For any $\epsilon$ put $L:= \log \epsilon$ and
set $\xi := e^{L_1} = e^N\epsilon$. Note that $\xi \rightarrow 0$ as
$\epsilon \rightarrow 0$. 

\begin{theorem}
\label{mainbound}
There is a bound $K$ depending on $u,w$ and a real constant $F$ (which will be $\zeta _0+2\nu -\zeta _1$ later on),
but with $K$ independent of $\epsilon$, $\xi$ and $\mu$,
with the following properties. 
Suppose $T$ is a furnished tree and $\Phi : T \rightarrow X$ is a flowing map
such that the top vertices are outside of any $D_w(P_i)$ and 
such that the bottom vertices are each mapped into some $D_{\epsilon}(P_i)$. 
Suppose furthermore that $g(z^{\rm bot} ) \geq g(z^{\rm top}) -F$. 
Then we can cut $T$ into a tree $T'$ onto which are attached subtrees,
such that $\Phi$ maps the bottom vertices of $T'$ into 
various $D_{\xi}(P_i)$ and such that the number of bottom vertices of $T'$ is
bounded by $K$. 
\end{theorem}
\begin{proof}
Among the subtrees that we strip off are any ones starting with edges
$e$ for which $i_e=j_e$. In particular we may assume from the start
that $T$ has no such edges. 

Next group the bottom vertices into series connected by intervals
where the bounding loop of the interval is mapped into $D_u(P)$.
There is a bound $K_1$ for the number of such series, because any loop 
which goes out of $D_u(P)$ has to contribute at least a  certain fixed
amount to $g(z^{\rm top}) -g(z^{\rm bot} )$. Next we can look at a specific
series. It is the set of bottom vertices of a subtree $T_1$ obtained
by taking the union of all of the loops joining the bottom vertices
together. Note in particular that $\Phi (T_1)\subset D_u(P)$.
Let $\kappa$ denote the boundary path of $T_1$. Note that the subtree $T_1$
doesn't necessarily include all strands emanating from all of its vertices.
However, if $v$ is a vertex on $\kappa$ corresponding to an adjoining edge $e$
not in $\kappa$, then either $e$ goes into the interior of the region 
bounded by $\kappa$, in which case $e$ starts a subtree mapped into $D_u(P)$
and such that all bottom edges go into $D_{\epsilon}(P)$; or else
it goes out of the region bounded by $\kappa$ in which case $e$ is not
a part of the tree $T_1$. In the former case, the normal form of Lemma
\ref{special} applies to the subtree starting at $e$. In the latter case,
the subtree starting at $e$ could be in normal form or not. However,
if $e$ is an edge going out of $\kappa$ such that the subtree starting
at $e$ is not in the normal form of Lemma \ref{special}, then this
subtree contains at least one strand which goes out of $D_u(P)$.
By Lemma \ref{outgoing} it also contains a strand which goes out of $D_w(P)$
and there is a global bound $K_2$ on the number of such 
edges $e$. If we cut $\kappa$ at vertices $v$ where such edges $e$ 
go out, then it is cut into $\leq K_2$ strands $\kappa '$ and each
little strand has only vertices corresponding to normal-form subtrees.
Finally, by the bound of Lemma \ref{innerbound} there are no more than $K_3$
such vertices on each little strand $\kappa '$. Each of these normal-form
subtrees can be cut at the point where it goes into $D_{\xi}(P)$,
and there is only one such point for each subtree. Thus if we trim off 
the tree $T_1$ at all of the points where the strands enter
$D_{\xi}(P)$, there are at most $K_2K_3$ bottom vertices. Finally,
since there were at most $K_1$ subtrees $T_1$ corresponding to series
of bottom vertices, we can trim off $T$ to a tree $T'$ where
there are at most $K_1K_2K_3$ bottom vertices, all going inside
some $D_{\xi}(P_i)$. 
This proves the theorem.
\end{proof}

\section{Proofs}
\label{proof}

By a {\em multisingular point} we mean a point $y=(y_1,\ldots , y_k)\in Z^{\ast}_I$ such that the $y_n$
are singular points of the functions $g_{i_ni_{n+1}}$. Note that in our situation the singular points in 
$Z$ are the preimages of the turning points $P\in \Rr$ corresponding to the places where the
Higgs field $\theta$ has singular eigenvalues. 

If $z$ is a point in $Z^{\ast}_I$ with ${\bf r}(z)\leq M$ then in particular each $z_i$ is at 
distance $\leq M$ from $p$ with respect to $d\sigma$ (using the notations of \S \ref{remoteness}). 
This defines a relatively compact subset of $Z$,
containing a finite number of singular points. It is improved with the lemma below. 

Define $S_M$ to be the set of complex values of the form $g(y)$ where $y$ are multisingular points
with ${\bf r}(y)\leq M_0+2M$. This is the subset which is to enter into the 
definition of analytic continuation with locally finite branching for $f(\zeta )$. 

\begin{lemma}
\label{finiteness}
For each $M$, the set $S_M$ is finite.
\end{lemma}
\begin{proof}
There is a positive constant $c$ such that if $P_1$ and $P_2$ are distinct singular points,
then the distance from $P_1$ to $P_2$ using the metric $d\sigma$ is at least $c$.
Suppose $y=(y_1, \ldots , y_k)$ is a multisingular point, so each $y_i$ is a singular point.
If ${\bf r}(y)\leq M_0+2M$ then, in view of the definition of ${\bf r}$ there are at most 
$(M_0+2M)/c$ indices $i$ such that $y_i\neq y_{i+1}$. Let $y':= (y'_1,\ldots , y'_{k'})$
be the sequence of distinct different points in the sequence $y$. Define a new multi-index
$I'$ by setting $i'_a=i_{b(a)}$ where $b(a)$ is the place with $y_{b(a)}=y'_a$ and
$y_{b(a)+1}=y'_{a+1}$. Then $y'\in Z_{I'}$ and $g_{I'}(y')=g_I(y)$. Now $k'\leq (M_0+2M)/c$ so
there are only a finite number of possibilities for $y'$ (the singular points themselves 
occuring in a fixed relatively compact subset of $Z$ as pointed out above). Thus there are 
only a finite number of possible values.  
\end{proof}

\begin{proof}[Proof of Theorem \ref{main}]

Suppose we have already analytically continued $f$
along a piecewise linear path of length $\leq M_1$. Inductively we may assume that the points 
of $\eta$ have remoteness $\leq M_0+ 2M_1$. If we add a segment $S$ then the total length of the path is
$\leq M$ where $M=M_1+|S|$. We assume that $S$ doesn't meet any of the points in $S_M$. 

Fix a number $\nu >0$ so that the segment $S$ stays at a distance $>2\nu $ away from the points of $S_M$.
Choose our neighborhoods $N_i$ with $N_1$ being the oval around $S$ of radius $\nu$, so $N_1$
stays at a distance $>\nu$ away from the points of $S_M$. 
Let $K$ be the bound of Theorem \ref{mainbound}. Choose $\epsilon$ small enough so that 
if $z\in D_{\xi}(P)$ then for any $i$ 
$$
\left| \int _P^z\alpha _i\right| < \frac{\nu}{K}.
$$
We show that all points of the chain $F\tau$ are sent (by $g$) outside of $N_3$. Suppose on the contrary that
we had a point, corresponding to a tree $T$, such that 
$g(z^{\rm bot}(T))\in N_3$. 

By Theorem \ref{mainbound} there exists a pruning $T'$ of $T$ with $\leq K$ bottom
vertices, such that
for every bottom vertex $v$ of $T'$ we have $\Phi (v)\in D_{\xi}(P(v))$ for some singular point $P(v)$.
In particular, the point $z^{\rm bot}(T')$ which is the vector of these $\Phi (v)$ is near to a point
$y=(\ldots , P(v),\ldots )$. More precisely we obtain from $k\leq K$ and the bound above,
$$
|g(y)-g(z^{\rm bot}(T'))| < \nu .
$$
On the other hand, if $g(z^{\rm bot}(T'))$ were inside $N_3$ then the singular point $y$ would occur
below points of $\eta$ at distance $\leq |S|$, and hence below points of $p$
at distance $\leq M_0+2M_1 + |S| < M$, therefore
$g(y)$ must be included in $S_M$. On the other hand, the point $g(z^{\rm bot}(T'))$ occurs
on the real segment between $g(z^{\rm bot}(T))$ and some point of $g(\eta ')$. This contradicts the assumption
that  the neighborhood $N_1$ stays away from $S_M$ by distance at least $\nu$. This shows that all points of
$g(F\tau )$ are outside of $N_3$, and completes the proof that we can analytically continue $f(\zeta )$ along
the segment $S$.  

Finally in order to maintain the inductive hypothesis we note that, cutting everything off
fairly close to the segment $S$ we can insure that the points of the new cycle of integration
$F\tau $ (and also $F\psi$) are remote from points of $\eta$ at distance $\leq 2|S|$, hence they have remoteness
$\leq M_0 + 2M$ as required. 
\end{proof}

\begin{proof}[Proof of Theorem \ref{mainplus}]---
The statement is essentially contained in that of Theorem \ref{main}, but we need to show that $S_M\subset \Sigma
(\gamma )$. In other words, if $z\in Z_I$ is a multisingular point, we need to show that $g(z)$ is the 
integral of the tautological form on a piecewise homotopy lifting $\tilde{\gamma}$. Recall the formula
$$
g_I(z)=\int _p^{z_1}\alpha _{i_0} + \ldots + \int _{z_k}^q \alpha _{i_k}.
$$
Let $\tilde{\gamma}'_i$ be the path joining $z_i$ to $z_{i+1}$ where by convention $z_0=p$ and $z_{k+1}=q$.
These paths are unique up to homotopy because we are working in the contractible universal cover $Z$.
Composing the main projection $Z\rightarrow Y$ with
Galois automorphisms of $Y$ and then the projection $Y\rightarrow V$, gives 
projections $\tau _i:Z\rightarrow V$ 
which commute with the projection to $X$, such that $\alpha _i$ is
the pullback of the tautological form $\alpha$ on $V$, i.e. $\alpha _i=\tau _i^{\ast}(\alpha )$.
We can put $\tilde{\gamma}_i := \tau _i\circ \tilde{\gamma}'_i$. The collection 
$\tilde{\gamma}= \{ \tilde{\gamma}_i \}$ is a piecewise homotopy lifting of $\gamma$. To see this, note
that the projections to $X$ of the $\tilde{\gamma}_i$ are equal to the projections of the original
$\tilde{\gamma}'_i$, so these join together to give a path homotopic to the projection of the path
from $p$ to $q$ in $Z$. Since the lifts $p,q\in Z$ were chosen to correspond to our original path $\gamma$ in
$X$, so the composite path in $X$ is homotopic to $\gamma$. Our formula for $g_I(z)$ becomes 
$$
g_I(z)=\sum _i \int _{\tilde{\gamma}_i}\alpha = \int _{\tilde{\gamma}}\alpha .
$$
This shows that $S_M$ is a subset of $\Sigma (\gamma )$.
\end{proof}

\section{Conclusion}
\label{conclusion}

We close with a few more general remarks about the consequences of Theorem \ref{main}. The first is
to note that it also applies to any polynomials in the transport matrix coefficients, in particular
to the Procesi coordinates for $M_B$. 

\begin{lemma}
\label{convolution}
Suppose $f_1$ and $f_2$ have analytic continuations with locally finite branching, then
the same is true for their convolution $f_1\ast f_2$. 
\end{lemma}
\begin{proof}
This was proven in \cite{CNP}.
See also the proof of \cite{abmade} Lemma 11.1.
There, the proof of locally finite branching for the convolutions uses only locally finite branching for
the two functions. 
\end{proof}

\begin{corollary}
\label{polynomials}
If $P(t)$ is a polynomial in the transport matrices for various paths, then 
then the Laplace transform of $P(t)$ has an analytic continuation with locally finite branching.
\end{corollary}
\begin{proof}
The Laplace transform of a 
product of functions $m_1(t)m_2(t)$ is the convolution of their Laplace transforms, so Lemma \ref{convolution}
and Theorem \ref{main} give the result.  
\end{proof}

The next remark is about 
the growth rate of $m(t)$. This is measured by the hull $\hull (m)$ defined in the introduction.

For reference we indicate first an elementary argument showing that $\hull (m)$
is convex. Indeed, if $\zeta _0$ is a point which is not in $\hull (m)$, then by definition there is an angular
sector ${\bf s}$ in which $m(t)e^{-\zeta _0t}$ is rapidly decreasing. Suppose $u$ is a complex number such
that $\zeta _0 + u$ is in $\hull (m)$. Again by the definition of $\hull (m)$ this implies that 
$m(t)e^{-\zeta _0t}e^{-ut}$ is no longer rapidly decreasing in 
any part of ${\bf s}$. This means that ${\bf s}$ is contained in the half-plane $\Re ut \leq 0$.
In particular, for any vector $u'$ which is a negative real multiple of $u$, we have that
$\Re u't \geq 0$ so  $m(t)e^{-\zeta _0t}e^{-u't}$ is rapidly decreasing on ${\bf s}$, therefore $\zeta _0+u'$
is not in $\hull (m)$. This  proves the convexity. 

Next we can characterize $\hull (m)$ as the intersection of all closed half-planes $H\subset \cc$ such that
the Laplace transform $f$ of $m$ admits an analytic continuation over the complementary open half-plane
(this would give another proof of convexity). 
Indeed, if a point $\zeta$ is in the complement of $\hull (m)$ then the sector along which 
$m(t)e^{-\zeta t}$ is rapidly decreasing provides an open half-plane containing $\zeta$ over which $f$ can
be analytically continued. This shows one inclusion. The other inclusion is clear from the inverse Laplace
transform.  

The hull is related to growth rates as follows. If $\hull (m)$ is a single point, then some
multiplicative translate of the form $m(t)e^{\zeta t}$ has sub-exponential growth. 
If $\hull (m)$ contains at least a line segment, then we say $m$ is {\em semistrictly exponential}: 
for sectors covering all but two directions we have a lower bound of the form $|m(t)| \geq ce^{a|t|}$,
and in particular there is a positive lower bound for the possible exponents $a$ which can enter into bounds of the form 
$|m(t)|\leq Ce^{a|t|}$. If $\hull (m)$ contains a nonempty interior then we say $m$ is {\em strictly exponential}:
there is a lower bound of the
form $|m(t)| \geq ce^{a|t|}$ valid in all directions. 

Unfortunately we are only able to show that some monodromy matrix is semistrictly exponential in the generic
case of Corollary \ref{hullcor} of the introduction. 

\begin{proof}[Proof of Corollary \ref{hullcor}]
By Theorem \ref{mainplus}, the Laplace transform has locally finite branching (Definition \ref{lfb}).
Choose $M$ big enough so that one goes all the way around $\hull (m)$ with a path of length $\leq M$.
Let $S^{\rm real}_M\subset S_M$ be the subset of non-removable singularities of the Laplace transform
attainable by a path of length $\leq M$
(which is finite because $S_M$ is finite).  
Then $f$ admits an analytic continuation to an open half-plane if and only if this half-plane doesn't meet
$S^{\rm real}_M$. Therefore $\hull (m)$ is a polygon. 

We show by specialization that for some fundamental group elements at least,
$\hull (m)$ is not reduced to a single point. General considerations using Hartogs' theorem
show that
if the monodromy is semistrictly exponential for a special curve going
to infinity, then 
the same will be true away from a piecewise holomorphic real codimension $2$ divisor.

We choose as special curve the family of connections on the 
trivial bundle of the form $d+B + tA$ with $A$ diagonal and $B$ off-diagonal,
everything being holmorphic on $X$, that was originally considered in \cite{abmade}.
In that case, we get asymptotic expansions whose coefficients can be calculated. 
One route is to note that for generic values of $A$ and $B$, calculation of the coefficients gives
nonzero coefficients at more than one singular point. Another route would be to note that if there
were only one singularity for the monodromy matrices for this family, then the monodromy representation would
actually have polynomial growth. That possibility is ruled out by specializing again to a direct sum of a $2\times 2$
system and trivial systems, and noting that for $2\times 2$ systems we have proven (in the paper \cite{tokyo}) that
the monodromy representation always has growth at least $e^{t^{1/k}}$ for some integer $k$. 

In any case by either of these two
routes we can conclude that the Laplace transform for at least one 
monodromy matrix has at least two singularities.
\end{proof}

It is perhaps more interesting to note that the same thing also works for the
Procesi coordinates.
This improves, at least for certain generic points at 
infinity approached from certain sectors, 
the bound given in \cite{tokyo}.

\begin{corollary}
\label{exponential}
For each family $(E,\nabla + t\theta )$ going to  infinity at a generic Higgs bundle $(E,\theta )$,
let $\rho _t$ denote the family of monodromy representations, thought of as a point in $M_B$.
Let $R_i:M_B\rightarrow \cc$ denote a set of  Procesi coordinates giving an affine embedding.
Write by abuse of notation $R_i(t):= R_i(\rho _t)$. Then each $\hull (R_i)$ is a polygon, and 
for general $(E,\theta )$ (in a dense open set)
at least one $R_i$ is semistrictly exponential (i.e. its hull has at least two vertices).
If we define $|\rho _t|:= \sup _i |R_i(t)|$ then for general $(E,\theta )$ and for
a family of sectors of $t\rightarrow \infty$
covering all but possibly two opposite directions, we have bounds of the form 
$$
|\rho _t| \geq ce^{a|t|}
$$
with $a>0$. 
\end{corollary}
\begin{proof}
The same proof as for Corollary \ref{hullcor} works here too.
\end{proof}

Lastly it is important to reiterate that, in spite of the above consequences, 
the result of Theorem \ref{main}
is highly unsatisfactory in that it doesn't say anything about the behavior of the Laplace transform $f(\zeta )$
near the singularities. It doesn't even seem clear what the answer will be: on the one hand one can 
imagine that an improvement of the present analysis, potentially based on Remark (i) following the proof of
Lemma \ref{modif}, might lead to a polynomial bound for the singularities. On the other hand, a crude look at
the present argument yields no such bound, and it is also quite concievable that the poles in the matrix $B$ 
lead unavoidably to more complicated singularities of $f(\zeta )$. This is undoubtedly true in the general case
where $B$ has poles of order $>1$. 

This problem also leads to the unsatisfactory statement of Corollary \ref{exponential}: if we could
calculate exactly where the singularities were we could probably show that for generic values of $(E,\theta )$ the
singularities would
span a convex hull with nonempty interior, in other words that the monodromy families $\rho _t$ would be
strictly exponential. This would be a more significant improvement of the result of \cite{tokyo}.

The result of Theorem \ref{main} should be thought of as a weak form of ``resurgence'' for the monodromy 
function $m(t)$ and its Laplace transform. The problem of 
getting more precise information about this behaviour is probably most naturally attacked using new ideas
and techniques for resummation such as have been developped by the school of J.-P. Ramis.



\begin{thebibliography}{55}

\bibitem{AokiKawaiTakei} 
T. Aoki, T. Kawai, Y. Takei. On the exact steepest descent method: A new method for the description of Stokes curves.
{\em J. Math. Phys.} {\bf 42} (2001), 3691-3713.

\bibitem{Arinkin}
D. Arinkin. Orthogonality of natural sheaves on moduli stacks of $SL(2)$-bundles with connections on $\pp ^1$ minus
$4$ points. {\em Selecta Math.} {\bf 7} (2001), 213-239. 

\bibitem{ArinkinLysenko}
D. Arinkin, S. Lysenko.
On the moduli of  ${\rm SL}(2)$-bundles with connections on $\pp ^1 -\{x_1,\ldots, x_4\}$. 
{\em Int. Math. Res. Not.} {\bf 19} (1997), 983-999.

\bibitem{BalianBloch}
R. Balian, C. Bloch.  Solution of the Schr\"{o}dinger equation in terms of classical paths. 
{\em Ann. Physics} {\bf 85} (1974), 514-545.

\bibitem{BenderEtAl}
C. Bender, M. Berry, P. Meisinger, V. Savage, M. Simsek. Complex WKB analysis of energy-level degeneracies
of non-Hermitian Hamiltonians. {\em J. Phys. A: Math. Gen.} {\bf 34} (2001), L31-L36. 

\bibitem{BerkNevinsRoberts}
H. Berk, W. Nevins, K. Roberts. New Stokes' line in WKB theory. {\em J. Math. Phys.} {\bf 23} (1982), 988-1002. 

\bibitem{Boutet}
L. Boutet de Monvel. {\em M\'ethodes R\'esurgentes: analyse alg\'ebrique des perturbations singuli\`eres}
{\sc Travaux en Cours}, Hermann, Paris (1994). 

\bibitem{RamisEtAl}
M. Canalis-Durand, J.-P. Ramis, R. Sch\"{a}fke, Y. Sibuya. Gevrey solutions of singularly
perturbed differential equations. {\em J. Reine Angew. Math.} {\bf 518} (2000), 95-129.

\bibitem{CNP}
B. Candelpergher,J.-C. Nosmas, F. Pham.
{\em Approche de la r\'esurgence}. {\sc Actualit\'es Math\'ematiques}, Hermann, Paris (1993). 

\bibitem{Chen}
K.-T. Chen. Integration of paths, geometric invariants and a generalized Baker-Hausdorff formula.
{\em Ann. of Math.} {\bf 65} (1957), 163-178.

\bibitem{CostinCostin}
O. Costin, R. Costin. Asymptotic properties of a family of solutions of the Painlevé equation $P_{VI}$. 
{\em Int. Math. Res. Notices} {\bf 22} (2002),  1167-1182.

\bibitem{CullerShalen}
M. Culler, P. Shalen. Varieties of group representations and splittings of $3$-manifolds. 
{\em Ann. of Math.} {\bf 117} (1983), 109-146.

\bibitem{DDthesis}
E. Delabaere, H. Dillinger. 
Contribution à la résurgence quantique. Résurgence de Voros et fonction spectrale de Jost.
Thesis, Universit\'e de Nice Sophia-Antipolis (1991). 


\bibitem{DDP}
E. Delabaere, H. Dillinger, F. Pham. R\'esurgence de Voros et p\'eriodes des courbes hyperelliptiques.
{\em Ann. Inst. Fourier} {\bf 43} (1993), 163-199. 

\bibitem{DDP2}
E. Delabaere, H. Dillinger, F. Pham. Exact semiclassical expansions for one-dimensional quantum oscillators.
{\em J. Math. Phys.} {\bf 38} (1997), 93-132. 

\bibitem{DelabaereHowls}
E. Delabaere, C. Howls. Global asymptotics for multiple integrals with boundaries. {\em Duke Math. J.} {\bf 112}
(2002), 199-264. 

\bibitem{DeligneLet} 
P. Deligne. Letter to the author.

\bibitem{Donagi}
R. Donagi.  Spectral covers. {\em Current topics in complex algebraic geometry (Berkeley, 1992/93)},
{\sc Math. Sci. Res. Inst. Publ.} {\bf 28}, Cambridge Univ. Press (1995), 65-86. 

\bibitem{EcalleOrsay}
J. Ecalle. Les fonctions r\'esurgentes. Tomes I, II, III. {\em Publications Math\'ematiques d'Orsay}
{\bf 81-5}, {\bf 81-6} (1981), {\bf 85-5} (1985). See also Orsay preprint 84T 62. 

\bibitem{Ecalle1}
J. Ecalle. The acceleration operators and their applications to differential equations, quasianalytic 
functions, and the constructive proof of Dulac's conjecture. 
{\em Proc. ICM-90, Kyoto}, vol. II, Springer (1991), 1249-1258.

\bibitem{EcalleAbordable}
J. Ecalle. Singularit\'es non abordables par la g\'eom\'etrie. {\em Ann. Inst. Fourier} {\bf 42} (1992), 73-164.

\bibitem{EcalleVallet}
J. Ecalle, B. Vallet. Correction and linearization of resonant vector fields and diffeomorphisms. {\em Math. Z.}
{\bf 229} (1998), 249-318.

\bibitem{GentileMastropietro}
G. Gentile, V. Mastropietro.  Methods for the analysis of the Lindstedt series for KAM tori and renormalizability in 
classical mechanics. A review with some applications.
{\em Rev. Math. Phys.} {\bf 8} (1996), 393-444. See also {\tt chao-dyn 9506004}.

\bibitem{Hain}
R. Hain. The de Rham homotopy theory of complex algebraic varieties, I. {\em $K$-theory} {\bf 1} (1987), 271-324.
 
\bibitem{Hitchin}
N. Hitchin. Stable bundles and integrable systems. {\em Duke Math. J.} {\bf 54}
(1987), 91-114.

\bibitem{Hitchin2}
N. Hitchin. The self-duality equations on a Riemann surface. {\em Proc. London Math. Soc.} {\bf 55}
(1987), 59-126. 

\bibitem{InabaIwasakiSaito}
Michi-aki Inaba, Katsunori Iwasaki, Masa-Hiko Saito.
Moduli of Stable Parabolic Connections, Riemann-Hilbert correspondence and Geometry of Painlev\'{e} equation of type VI, Part I.
Preprint {\tt math.AG/0309342}.

\bibitem{Kamvissis}
S. Kamvissis. Desingularization of a hyperelliptic curve associated with a doubly periodic Dirac potential. 
{\em Bull. Greek Math. Soc.} {\bf 46} (2002), 141-145.

\bibitem{Kanev}
V. Kanev. Spectral curves, simple Lie algebras, and Prym-Tjurin varieties. 
{\em Bowdoin 1987}, {\sc AMS Proc. Sympos. Pure Math.} {\bf 49}, Part 1 (1989), 627-645.

\bibitem{LosevManin}
A. Losev, Y. Manin. New moduli spaces of pointed curves and pencils of flat connections. 
{\em Michigan Math. J.} {\bf 48} (2000), 443-472.

\bibitem{MalgrangeRamis}
B. Malgrange, J.-P. Ramis. Fonctions multisommables. {\em Ann. Inst. Fourier} {\bf 42} (1992), 353-368. 

\bibitem{Procesi}
C.  Procesi. The invariant theory of $n\times n$ matrices. 
{\em Advances in Math.} {\bf 19} (1976), 306-381.

\bibitem{SanguinettiWoodhouse}
G. Sanguinetti, N. Woodhouse. Geometry of dual isomonodromic deformations. 
\newline
Preprint,
\verb}http://www.maths.ox.ac.uk/~nwoodh/dual.pdf}

\bibitem{abmade}
C. Simpson. {\em Asymptotic Behavior of Monodromy}, {\bf L. N. M. 1502}, Springer (1991).

\bibitem{hbls}
C. Simpson. Higgs bundles and local systems. {\em Publ. Math. I.H.E.S.} {\bf 75} (1992), 5-95. 

\bibitem{moduli}
C. Simpson. Moduli of representations of the fundamental group of a smooth projective variety, I, II.
{\em Publ. Math. I.H.E.S.} {\bf 79} (1994), 47-129 and {\bf 80} (1994), 5-79. 

\bibitem{naht}
C. Simpson. Nonabelian Hodge theory. {\em Proceedings ICM-90, Kyoto}, Springer, Tokyo (1991), 198-230. 

\bibitem{tokyo}
C. Simpson. A lower bound for the size of monodromy of systems of ordinary differential equations.
{\em Algebraic Geometry and Analytic Geometry---ICM-90 Satellite Conference Proceedings}, Springer
(1991), 198-230.

\bibitem{hfnac}
C. Simpson. The Hodge filtration on nonabelian cohomology. {\em AMS Proceedings of Symposia in Pure Mathematics}
{\bf 62.2} (1997), 217-281.

\bibitem{Stenger}
C. Stenger. Points tournants de syst\`emes d'\'equations diff\'erentielles ordinaires singuli\`erement perturb\'ees.
Thesis, Universit\'e Louis Pasteur, Strasbourg (1999). 
\newline
{\tt http://www-irma.u-strasbg.fr/irma/publications/1999/99019.shtml}

\bibitem{Voros}
A. Voros. The return of the quartic oscillator: the complex WKB method. {\em Ann. Inst. Henri Poincar\'e, Physique Th\'eorique}
{\bf 39} (1983), 211-338. 

\bibitem{VorosBourbaki}
A. Voros. Probl\`eme spectral de Sturm-Liouville: le cas de l'oscillateur quartique. {\em S\'eminaire Bourbaki}
{\bf 602} (1982/83), 95-104.

\bibitem{Voros2}
A. Voros. R\'esurgence quantique. {\em Annales de l'institut Fourier} {\bf 43} (1993), 1509-1534.

\bibitem{Woodhouse}
N. Woodhouse. The symplectic and twistor geometry of the general isomonodromic deformation problem. 
{\em J. Geom. Phys.} {\bf 39} (2001), 97-128. 
\newline
See also the slides at 
\verb}http://www.maths.ox.ac.uk/~nwoodh/ini.pdf}

\bibitem{ZinnJustin}
J. Zinn-Justin. Analyse des instantons et r\'esultats exacts. {\em Ann. Inst. Fourier} 
{\bf 53} (2003), 1259-1285 (to appear). 


\end{thebibliography}
\end{document}